\documentclass[11pt]{article}
\usepackage{amsmath,amsfonts}
\usepackage{amsthm,amssymb}
\usepackage{epsfig}
\usepackage{leftidx}
\usepackage{graphicx}
\usepackage{enumitem}
\usepackage{color, epstopdf}
\usepackage{cite}
\usepackage{graphicx,wrapfig,subfigure}

\newtheorem{theorem}{Theorem}[section]
\newtheorem{lemma}{Lemma}[section]

\newtheorem{example}{Example}
\newtheorem{remark}{Remark}

\newcommand{\dfd}[1]{\partial^{#1}_\tau}

\newcommand{\diff}{\triangledown_\tau}
\newcommand{\piff}{\partial_\tau}
\newcommand{\defeq}{:=}
\newcommand{\zd}{\,\mathrm{d}}
\newcommand{\abs}[1]{\left|#1\right|}

\newcommand{\bra}[1]{\left(#1\right)}
\newcommand{\brab}[1]{\big(#1\big)}
\newcommand{\braB}[1]{\Big(#1\Big)}
\newcommand{\brat}[1]{(#1)}
\newcommand{\kbra}[1]{\left[#1\right]}
\newcommand{\kbrab}[1]{\big[#1\big]}
\newcommand{\kbraB}[1]{\Big[#1\Big]}
\newcommand{\kbrat}[1]{[#1]}
\newcommand{\myinner}[1]{\left\langle#1\right\rangle}
\newcommand{\myinnerb}[1]{\big\langle#1\big\rangle}

\newcommand{\mynorm}[1]{\left\|#1\right\|}
\newcommand{\mynormb}[1]{\big\|#1\big\|}

\title{Asymptotically compatible energy of variable-step fractional 
	BDF2 formula for time-fractional Cahn-Hilliard model}
\author{Hong-lin Liao\thanks{ORCID 0000-0003-0777-6832. College of Mathematics, 
		Nanjing University of Aeronautics and Astronautics,
		Nanjing 211106, China; Key Laboratory of Mathematical Modeling
		and High Performance Computing of Air Vehicles (NUAA), MIIT, Nanjing 211106, China. 
		Emails: liaohl@nuaa.edu.cn and liaohl@csrc.ac.cn.
		This author's work is supported by NSF of China
		under grant number 12071216.}
	\and
	Nan Liu\thanks{College of Mathematics, Nanjing University of Aeronautics and Astronautics,
		Nanjing 211106, China. Email: liunan@nuaa.edu.cn.}		
	\and 
	Xuan Zhao\thanks{School of Mathematics, Southeast University, Nanjing 210096, China. 
		Email: xuanzhao11@seu.edu.cn.}
}
\begin{document}

\maketitle

\begin{abstract}
	A new discrete energy dissipation law of the variable-step
	fractional BDF2 (second-order backward differentiation formula) 
	scheme is established for time-fractional Cahn-Hilliard model
	with the Caputo's  fractional derivative of order $\alpha\in(0,1)$,
	under a weak step-ratio constraint $0.4753\le \tau_k/\tau_{k-1}<r^*(\alpha)$,
	where $\tau_k$ is the $k$-th time-step size and $r^*(\alpha)\ge4.660$ for $\alpha\in(0,1)$.
	We propose a novel discrete gradient structure by a local-nonlocal splitting technique, 
	that is, the fractional BDF2 formula is split into 
	a local part analogue to the two-step backward differentiation formula of the first derivative 
	and a nonlocal part analogue to the L1-type formula of the Caputo's  derivative.
	More interestingly, in the sense of the limit $\alpha\rightarrow1^-$,  
	the discrete energy and the corresponding energy dissipation law are asymptotically compatible with the associated discrete energy and the energy dissipation law of the variable-step BDF2 method for the classical Cahn-Hilliard equation, respectively. Numerical examples with an adaptive stepping procedure are provided to demonstrate the accuracy and the effectiveness of our proposed method.
\noindent{\bf AMS subject classiffications.}\;\; 35Q99, 65M06, 65M12, 74A50
\end{abstract}

\section{Introduction}
\setcounter{equation}{0}

Linear and nonlinear diffusion equations with fractional time derivatives
have become widely models describing anomalous
diffusion processes \cite{Hilfer:2000,KilbasSrivastavaTrujillo:2006,TangYuZhou:2019}. 
These models always exhibit multi-scaling time behaviour, 
which makes them suitable for the description of different diffusive regimes 
and characteristic crossover dynamics in complex systems.
To capture the multi-scale behaviors
in time-fractional differential equations,
adaptive time-stepping strategies, namely, small time steps are utilized when the solution varies
rapidly and large time steps are employed otherwise, would be practically useful especially in 
the simulations of long-time coarsening dynamics 
\cite{JiLiaoGongZhang:2019,JiZhuLiao:2022TFCH,LiaoJiWangZhang:2021ch,
	LiaoZhang:2021MC,LiaoTangZhou:2020Anenergy,LiaoZhuWang:2021}. 
It is natural to require practically and theoretically reliable
time-stepping methods on general setting of time step-size variations, 
or on a wider class of nonuniform time meshes.

This work develops a new discrete energy dissipation law of the second-order
fractional backward differentiation formula (FBDF2) for the 
time-fractional Cahn-Hilliard (TFCH) flow in modelling 
the coarsening process with a general coarsening rate \cite{AlMaskariKaraa:2021,FritzRajendranWohlmuth:2022,
	LiuChengWangZhao:2018,QuanWang:2022jcp,QuanTangWangYang:2022,TangYuZhou:2019,
	JiZhuLiao:2022TFCH,TangWangYang:2022},
\begin{align}\label{cont:TFCH model}
	\partial_t^\alpha \Phi=\kappa\Delta\mu\quad\text{with}\quad
	\mu:=\tfrac{\delta E}{\delta \Phi}=f(\Phi)-\epsilon^2\Delta\Phi,
\end{align}
where $f(\Phi):=F'(\Phi)$ and $E[\Phi]$ is the Ginzburg-Landau energy functional
\begin{align}\label{cont:classical free energy}
	E[\Phi]:=\int_{\Omega}\braB{\frac{\epsilon^2}{2}\abs{\nabla\Phi}^2
		+F(\Phi)}\zd\mathbf{x}\quad\text{with the potential $F(\Phi):=\frac{1}{4}\bra{\Phi^2-1}^2$}.
\end{align}
The notation $\partial_t^\alpha:={}_{0}^{C}\!D_{t}^{\alpha}$ represents
the Caputo's derivative of order $\alpha\in(0,1)$, defined by
\begin{align}\label{def:Caputo derivative}
	(\partial_{t}^{\alpha}v)(t)
	:=\int_{0}^{t}\omega_{1-\alpha}(t-s)v'(s)\zd{s}
	\quad \text{with $\omega_{\beta}(t):=t^{\beta-1}/\Gamma(\beta)$ for $\beta>0$.}
\end{align} 
The real valued function $\Phi$ represents the concentration
difference in a binary system on the domain $\Omega\subseteq\mathbb{R}^2$,
$\epsilon>0$ is an interface width parameter, and  $\kappa>0$ is the mobility coefficient.

Let $\myinner{\cdot,\cdot}$ and $\mynorm{\cdot}$ be the
$L^2(\Omega)$ inner product and the associated norm, respectively.
Always we use the standard norms of the Sobolev space $H^{m}\bra{\Omega}$
and the $L^p\bra{\Omega}$ space.
For $v$ and $w$ belonging to the zero-mean space
$\mathbb{\mathring V}:=\big\{v\in L^2\bra{\Omega}\,|\, \myinner{v,1}=0\big\}$,
we use the $H^{-1}$-like inner product $\myinner{v,w}_{-1}:=\myinner{\brat{-\Delta}^{-1}v,w}$ and
the induced norm $\mynorm{v}_{-1}:=\sqrt{\myinner{v,v}_{-1}}$. 
The well-possness of the TFCH model 
\eqref{cont:TFCH model} has been established in \cite{AlMaskariKaraa:2021,FritzRajendranWohlmuth:2022}
and the analysis indicated that the solution lacks the smoothness near the initial 
time while it would be smooth away from $t=0$,
also see \cite{DuYangZhou:2020,FritzRajendranWohlmuth:2022}.
The motions of interfaces and coarsening rates governed
by the TFCH equation were studied by Tang, Wang and Yang \cite{TangWangYang:2022}.
For the constant mobility, the sharp interface limit model is a fractional Stefan problem at 
the timescale $t=O(1)$; while on the timescale $t=O(\epsilon^{-1/\alpha})$, the sharp interface limit model is a fractional Mullins–Sekerka model. The TFCH model
with constant mobility was proven to preserve an $O(t^{-\alpha/3})$ coarsening rate, which is asymptotically compatible (as $\alpha\rightarrow1^-$)
with the coarsening rate $O(t^{-1/3})$ of the following classical Cahn–Hilliard (CH) equation,
\begin{align}\label{cont: CH model}
	\partial_t \Phi=\kappa\Delta\mu\qquad \text{for $t>0$}.
\end{align}  

As is well known, the CH model \eqref{cont: CH model} conserves the initial volume $\myinner{\Phi(t),1}=\myinner{\Phi(0),1}$ and
has the following energy dissipation law
\begin{align}\label{cont: CH energy decay law}
\frac{\zd E}{\zd t}+\kappa\mynorm{\nabla\mu}^2=0\qquad
\text{for $t>0$}.
\end{align}
Some efforts have been made to seek the reasonable versions of the energy $E[\Phi]$
and the energy dissipation law  \eqref{cont: CH energy decay law} 
for the TFCH equation \eqref{cont:TFCH model}. 
Tang, Yu and Zhou \cite{TangYuZhou:2019} established a global energy dissipation law, 
$E\kbrat{\Phi(t)}\le E\kbrat{\Phi(0)}$ for $t>0$, see also \cite{JiLiaoGongZhang:2019,QuanWang:2022jcp}.
%\begin{align*}
%	E\kbrat{\Phi(t)}\le E\kbrat{\Phi(0)} \quad \text{for $t>0$}.
%\end{align*}
In \cite{QuanTangYang: 2020csiam}, Quan, Tang and Yang derived a nonlocal energy decaying law, $\partial_{t}^{\alpha}E\le0$,
and a weighted energy dissipation law, $\partial_{t}E_{\eta}\le 0$, 
for a nonlocal energy $E_{\eta}(t):=\int_0^1\eta(\theta)E(\theta t)\zd \theta$, also cf.  
the dissipation-preserving augmented energy in \cite{FritzRajendranWohlmuth:2022}. 
They are quite different from the classical energy law
\eqref{cont: CH energy decay law}.
By reformulating  the Caputo form \eqref{cont:TFCH model} into the Riemann-Liouville 
form and applying the equality 
in \cite[Lemma 1]{AlsaediAhmadKirane:2015}, one can obtain the following
 variational energy dissipation law, cf. \cite[Section 1]{LiaoTangZhou:2020Anenergy},
\begin{align}\label{cont:TFCH variational energy law}
	\frac{\zd E_\alpha}{\zd t}
	+\frac{\kappa}{2}\omega_{\alpha}(t)\mynorm{\nabla\mu}^2
	-\frac{\kappa}{2}\int_0^t\omega_{\alpha-1}(t-s)\mynorm{\nabla\mu(t)-\nabla\mu(s)}^2\zd{s}=0,
\end{align}
where the variational (modified) energy 
\begin{align}\label{cont:TFCH variational energy}
E_\alpha[\Phi]:=E[\Phi]
	+\frac{\kappa}{2}\mathcal{I}_{t}^{\alpha}\mynorm{\nabla\mu}^2
	\quad \text{for $t>0$}.
\end{align}
We remark that the variational energy law \eqref{cont:TFCH variational energy law} is naturally consistent 
with the energy dissipation 
law \eqref{cont: CH energy decay law} of the CH model. Actually, the discrete energy dissipation laws 
of the variable-step L1 or L1$_R$ time-stepping schemes in 
\cite{JiZhuLiao:2022TFCH,LiaoTangZhou:2020Anenergy,	LiaoZhuWang:2021} %ZhuLiao:2022
were proven to be asymptotically compatible (in the limit $\alpha\rightarrow1^-$)
with the associated discrete energy laws of their integer-order counterparts.
However, the variational energy $E_\alpha[\Phi]$ is not asymptotically compatible with the
Ginzburg-Landau energy $E[\Phi]$.
Very recently, Quan \emph{et al.} \cite{QuanTangWangYang:2022} 
applied an equivalent definition (obtained
from the integration by parts) of the Caputo derivative  to find a new energy dissipation law (in our notations)
 \begin{align}\label{cont:TFCH QuanTangWangYang energy law}
 	\frac{\zd\widetilde{E}_\alpha}{\zd t}-\frac{\omega_{-\alpha}(t)}{2\kappa}\mynorm{\Phi(t)-\Phi(0)}_{-1}^2
 	+\frac{1}{2\kappa}\int_0^t\omega_{-\alpha-1}(t-s)\mynorm{\Phi(t)-\Phi(s)}_{-1}^2\zd{s}=0,
 \end{align}
where the modified energy $\widetilde{E}_\alpha[\Phi]$ is defined by
\begin{align}\label{cont:TFCH QuanTangWangYang energy}
	\widetilde{E}_\alpha[\Phi]:=E[\Phi]
	+\frac{\omega_{1-\alpha}(t)}{2\kappa}\mynorm{\Phi(t)-\Phi(0)}_{-1}^2
	-\frac{1}{2\kappa}\int_0^t\omega_{-\alpha}(t-s)\mynorm{\Phi(t)-\Phi(s)}_{-1}^2\zd{s}.
\end{align}
An interesting property is that the new modified energy $\widetilde{E}_\alpha[\Phi]$ is asymptotically compatible with the Ginzburg-Landau energy \cite[Proposition 5.1]{QuanTangWangYang:2022}, that is, $\widetilde{E}_\alpha[\Phi]\rightarrow E[\Phi]$ as $\alpha\rightarrow1^-$. They also derived the modified energy dissipation laws at discrete time levels for 
the L1 and L2-type implicit-explicit stabilized schemes on the uniform time mesh. Nonetheless, 
it is noticed that the modified discrete energies of these time-stepping schemes are
 not asymptotically compatible (except the steady state case) 
 with the associated discrete energies of their integer-order counterparts, 
 see more details in \cite[Theorem 5.2]{QuanTangWangYang:2022} or \cite[Theorems 4.1 and 4.2]{QuanTangWangYang:2022}.
% Moreover, it seems that the modified energy dissipation law \eqref{cont:TFCH QuanTangWangYang energy law}
% is not asymptotically compatible  with the energy dissipation law \eqref{cont: CH energy decay law}.

The above mentioned variable-step schemes are built by the piecewise linear 
or piecewise constant interpolating polynomial 
and only have a low order accuracy of $O(\tau^{2-\alpha})$ or $O(\tau^{1+\alpha})$ in time, see the related error analysis in \cite{Kopteva:2019,LiaoLiZhang:2018,LiaoMcLeanZhang:2019,Mustapha:2011,MustaphaAlMutawa:2012}.
In this paper, we shall build a {\em discrete gradient structure} (DGS) of fractional BDF2 formula
on general nonuniform meshes and establish an asymptotically compatible discrete energy dissipation law  at each time level  with an asymptotically compatible energy. 
As far as we know, there are few studies on the high-order energy stable schemes with unequal time-step sizes for time-fractional phase field models, see \cite{JiLiaoGongZhang:2019,QuanWang:2022jcp,QuanWu:2022}.

Application of quadratic interpolating polynomial to obtain the high-order approximations of Caputo derivative \eqref{def:Caputo derivative} was suggested and investigated in recent years, see \cite{GaoSunZhang:2014,LiaoMcLeanZhang:2021alikhanov,LvXu:2016}. The so-called L2-type (fractional BDF2) methods in \cite{GaoSunZhang:2014,LvXu:2016} were shown to be 
$(3-\alpha)$-order accurate on the uniform mesh for sufficiently smooth solutions. 
By using a nonuniform mesh such as the well-known graded meshes \cite{Kopteva:2019,Kopteva:2021,LiaoLiZhang:2018,LiaoMcLeanZhang:2019,
	LiaoMcLeanZhang:2021alikhanov,QuanWu:2022},
we will restore the optimal convergence when the solution is not smooth near the initial time. 
%This idea was tested in \cite{LiaoLyuVongZhao:2017} to resolve the
%initial singularity for the subdiffusion problem. 
Recently, Kopteva \cite{Kopteva:2021} applied the
inverse-monotonicity of discrete fractional-derivative operator to obtain  
sharp pointwise-in-time error bounds on quasi-graded meshes
for the linear subdiffusion equation and derived a mild constraint of 
grading parameter to recover the optimal order convergence
at any positive time.

We study the variable-step fractional BDF2 scheme on a general class of time meshes.
Consider the nonuniform time levels $0=t_{0}<t_{1}<\cdots<t_{k-1}<t_{k}<\cdots<t_{N}=T$ with the time-step sizes $\tau_{k}:=t_{k}-t_{k-1}$ for $1\leq{k}\leq{N}$ and the maximum time-step size 
$\tau:=\max_{1\leq{k}\leq{N}}\tau_{k}$. Also, let  $t_{k-1/2}:=(t_k+t_{k-1})/2$ for $k\ge1$, 
and the adjacent time-step ratio $r_1:=0$ and $r_k:=\tau_k/\tau_{k-1}$ for $k\ge2$.
%$$r_1:=0\quad\text{and}\quad r_k:=\tau_k/\tau_{k-1}\quad\text{for $k\ge2$}.$$
For any time sequence $v^k=v(t_k)$, define the backward
difference~$\diff v^k:=v^k-v^{k-1}$ and the difference  quotient~$\piff v^k:=\diff v^k/\tau_k$.
%Let $\Pi_{1,k}v$ denote the linear interpolant of a function~$v$ with
%respect to the nodes $t_{k-1}$~and $t_k$, and let $\Pi_{2,k}v$ denote the
%quadratic interpolant with respect to three nodes $t_{k-1}$, $t_k$~and $t_{k+1}$ for $k\ge1$.
%Always, let $\widetilde{\Pi_{\nu,k}}v:=v-\Pi_{\nu,k}v$  ($\nu=1,2$) 
%be the corresponding interpolation errors.
Let $\Pi_{2,k}v$, $k\ge1$, denote the
quadratic interpolant with respect to three nodes $t_{k-1}$, $t_k$~and $t_{k+1}$ for $k\ge1$.
It is easy to find (for instance, by using
the Newton forms of the interpolating polynomials) that 
\begin{align*}
	\bra{\Pi_{2,k}v}'=\piff v^k
	+\frac{2(t-t_{k-1/2})}{\tau_k+\tau_{k+1}}
	\brab{\piff v^{k+1}-\piff v^k}
	=\piff v^{k+1}+\frac{2(t-t_{k+1/2})}{\tau_k+\tau_{k+1}}
	\brab{\piff v^{k+1}-\piff v^k}.
\end{align*}
For any time-level $t_n$ with $n\ge2$, applying the quadratic interpolating polynomial $\Pi_{2,k}v$,
we have the following $(3-\alpha)$-order BDF2 type formula 
\cite{GaoSunZhang:2014,LvXu:2016} of Caputo derivative \eqref{def:Caputo derivative},
\begin{align}\label{def: FBDF2 formula}
	(\dfd{\alpha}v)^{n}
	\defeq&\,
	\int_{t_{n-1}}^{t_n}\omega_{1-\alpha}(t_{n}-s)\bra{\Pi_{2,n-1}v}'(s)\zd{s}
	+\sum_{k=1}^{n-1}\int_{t_{k-1}}^{t_k}\omega_{1-\alpha}(t_{n}-s)\bra{\Pi_{2,k}v}'(s)\zd{s}\nonumber\\
	\triangleq&\, \sum_{k=1}^{n}a_{n-k}^{(n)}\diff v^k
	+\frac{r_n\eta^{(n)}_{0}}{1+r_{n}}\brab{\diff v^{n}-r_n\diff v^{n-1}}
	+\sum_{k=1}^{n-1}\frac{\eta^{(n)}_{n-k}\brab{\diff v^{k+1}-r_{k+1}\diff v^k}}{r_{k+1}(1+r_{k+1})}
\end{align}
for $n\ge2$, where the positive coefficients $a_{n-k}^{(n)}$ and $\eta_{n-k}^{(n)}$ are defined by
\begin{align}
	&a^{(n)}_{n-k}\defeq\frac{1}{\tau_k}\int_{t_{k-1}}^{t_k}\omega_{1-\alpha}(t_{n}-s)\zd{s}
	\quad \text{for $1\le k\leq n$,}\label{coeff: ank}\\
	&\eta^{(n)}_{n-k}\defeq\frac{2}{\tau_{k}}\int_{t_{k-1}}^{t_k}
	\frac{s-t_{k-1/2}}{\tau_k}\omega_{1-\alpha}(t_{n}-s)\zd{s}\quad \text{for $1\le k\leq n$}.	
	\label{coeff: eta nk}
\end{align}
To start with the stepping formula \eqref{def: FBDF2 formula}, we use the standard L1 formula 
at the first time level,
\begin{align}\label{def: FBDF2 formula-start}
	(\dfd{\alpha}v)^{1}
	\triangleq&\, a_{0}^{(1)}\diff v^1\quad \text{with}\quad
	a^{(1)}_{0}\defeq\frac{1}{\tau_1}\int_{t_{0}}^{t_1}\omega_{1-\alpha}(t_{1}-s)\zd{s}.
\end{align}

Notice that if the fractional order $\alpha\to1^-$, then
$\omega_{3-\alpha}(t)\to t$, $\omega_{2-\alpha}(t)\to1$
and $\omega_{1-\alpha}(t)\to0$, uniformly for~$t>0$.  Thus,
$a^{(n)}_0=\omega_{2-\alpha}(\tau_n)/\tau_n\to1/\tau_n$ and 
$\eta^{(n)}_{0}=\alpha\tau_n^{-\alpha}/\Gamma(3-\alpha)\to 1/\tau_{n}$,
%$$a^{(n)}_0=\frac{\tau_n^{-\alpha}}{\Gamma(2-\alpha)}\to\frac1{\tau_n}\quad \text{and}\quad 
%\eta^{(n)}_{0}=\frac{\alpha\tau_n^{-\alpha}}{\Gamma(3-\alpha)}\to \frac1{\tau_n},$$
whereas $a^{(n)}_{n-k}\to0$ and $\eta^{(n)}_{n-k}\to0$ for $1\le k\le n-1$.
We shall call \eqref{def: FBDF2 formula} together with \eqref{def: FBDF2 formula-start}  as 
the fractional BDF2 (FBDF2) formula for brevity,
since the approximation \eqref{def: FBDF2 formula} degrades into 
the standard BDF2 formula 
\cite{LiLiao:2022,LiaoJiWangZhang:2021ch,LiaoTangZhou:2020bdf2,LiaoZhang:2021MC} 
of the derivative $\partial_t v$, that is, %$(\dfd{\alpha}v)^{n}\to D_2v^n$ as $\alpha\to1^-$.
\begin{align}\label{eq: BDF2 formula}
	(\dfd{\alpha}v)^{n}\to D_2v^n
	:=\frac{1+2r_n}{(1+r_{n})\tau_n}\diff v^n-\frac{r_n^2}{(1+r_{n})\tau_n}\diff v^{n-1}
	\qquad\text{as $\alpha\to1^-$.}
\end{align}

We can reformulate \eqref{def: FBDF2 formula} and \eqref{def: FBDF2 formula-start} into a compact form, 
\begin{align}\label{def: FBDF2 formula-compact}
	 (\dfd{\alpha}v)^{n}\triangleq\sum_{k=1}^{n}B_{n-k}^{(n)}\diff v^k\quad\text{for $n\ge1$,}
\end{align}
where the kernels $B_{n-k}^{(n)}$ are determined via \eqref{def: FBDF2 formula} and \eqref{def: FBDF2 formula-start}.
The above asymptotic property \eqref{eq: BDF2 formula} arises the main difficulty in the numerical analysis 
of the FBDF2 formula \eqref{def: FBDF2 formula}. That is, the second kernel $B_1^{(n)}$ 
would be negative when $\alpha$ is close to $1$, and the discrete kernels $B_{n-k}^{(n)}$ lose 
the monotonicity. The established stability and convergence theory \cite{Kopteva:2019,LiaoLiZhang:2018,LiaoMcLeanZhang:2019,LiaoMcLeanZhang:2021alikhanov}
for the variable-step L1 and L2-1$_{\sigma}$ formulas can not be applied here directly.  
Recently, the uniform FBDF2 formula was showed in \cite{QuanWang:2022jcp} 
to be positive semidefinite with the BDF2 multiplier $D_2v^n$, that is, $\sum_{k=1}^n\brat{D_2v^k}(\dfd{\alpha}v)^{k}\ge0$.
%\begin{align*}
%	\sum_{k=1}^n\brat{D_2v^k}(\dfd{\alpha}v)^{k}\ge0.
%\end{align*}
By combining with the recent scalar auxiliary variable technique, 
the resulting time-stepping scheme was proven to preserve a global energy law \cite{TangYuZhou:2019} of 
time-fractional phase-field equations. Quan and Wu \cite{QuanWu:2022} investigated 
the $H^1$ stability of the FBDF2 formula \eqref{def: FBDF2 formula} 
on general nonuniform time meshes for a linear subdiffusion equation and established
the positive semidefiniteness with the difference multiplier $\diff v^n$, that is,
$\sum_{k=1}^n\brat{\diff v^k}(\dfd{\alpha}v)^{k}\ge0$,
%\begin{align*}
%	\sum_{k=1}^n\brat{\diff v^k}(\dfd{\alpha}v)^{k}\ge0
%\end{align*}
under the following step-ratio condition
\begin{align}\label{eq: QuanWu Ratio condition FBDF2}
	0.4573328\le r_k\le 3.5615528\quad\text{for $k\ge2$.}
\end{align}
Nonetheless, these results in \cite{QuanWang:2022jcp,QuanWu:2022} would be inadequate to build
an asymptotically compatible discrete energy law for the TFCH equation \eqref{cont:TFCH model}.

The current work is inspired by the following DGS
\cite{LiaoJiWangZhang:2021ch} of the BDF2 formula,
\begin{align}\label{ieq: DGS-BDF2 formula}
	\brat{\diff v^n} D_2v^n
	\ge\frac{r_{n+1}^{3/2}}{2(1+r_{n+1})\tau_n}\brat{\diff v^n}^2
	-\frac{r_{n}^{3/2}}{2(1+r_{n})\tau_{n-1}}\brat{\diff v^{n-1}}^2\quad\text{for $n\ge2$.}
\end{align}
This DGS holds if the adjacent step-ratio restriction $0<r_k<r^\star\approx 4.864$ for $k\ge2$,
where $r^\star$ is the positive root of $1+2r^\star-(r^\star)^{3/2}=0$.
This DGS has been proven 
to be useful to establish the discrete energy dissipation laws at each time level of
the BDF2 scheme for the classical gradient flows \cite{LiaoJiWangZhang:2021ch}.
 It is natural to ask {\em whether the FBDF2 formula \eqref{def: FBDF2 formula} 
has a similar gradient structure under some constraint on the step-ratio $r_k$, 
and whether the resulting FBDF2 time-stepping for time-fractional gradient flows 
also has an energy dissipation law at discrete time levels.}

In the next section, we update the step-ratio stability constraint 
\eqref{eq: QuanWu Ratio condition FBDF2} into
\begin{align*}
	0.4753\le r_k< r^*(\alpha) \quad\text{for $k\ge2$ with $r^*(\alpha)\ge4.660$,}
\end{align*}
 and establish a novel DGS for the FBDF2 formula \eqref{def: FBDF2 formula}.
%\begin{align*}
%	(\diff v^n)(\dfd{\alpha}v)^{n}\geq \mathcal{G}\kbrab{\diff v^n}-\mathcal{G}\kbrab{\diff v^{n-1}}
%	+\mathcal{G}_R\kbrab{\diff v^n}
%	\quad \text{for $n\ge2$.}
%\end{align*}
Our main tools include a novel local-nonlocal splitting technique and 
the step-scaling technique so that the present analysis is quite different and 
would be more concise than those in \cite{QuanTangWangYang:2022,QuanWang:2022jcp,QuanWu:2022}. 
As an application, we consider an implicit FBDF2 stepping 
for the TFCH model \eqref{cont:TFCH model} in Section 3 and
 prove that it is uniquely solvable and has
 an asymptotically compatible energy law  at each time level. 
 More importantly, our discrete  energy is also asymptotically compatible with 
 the associated discrete energy of the BDF2 scheme for the CH equation.
 Numerical tests show that the proposed method is accurate with 
 an order of $O(\tau^{3-\alpha})$ in time. An adaptive time-stepping 
 procedure is also provided to speedup the simulations of long-time coarsening dynamics.

\section{Discrete gradient structure of FBDF2 formula}
\setcounter{equation}{0}
   
To derive the DGS of the nonuniform FBDF2 formula  \eqref{def: FBDF2 formula}, 
we propose a local-nonlocal splitting technique by splitting it into two parts
\begin{align}\label{eq: FBDF2-local nonlocal splitting}
	(\dfd{\alpha}v)^{n}
	\triangleq&\,\underbrace{\braB{\frac{1}{2-\alpha}a_{0}^{(n)}
			+\frac{r_n\eta^{(n)}_{0}}{1+r_{n}}}\diff v^{n}-\frac{r_n^2\eta^{(n)}_{0}}{1+r_{n}}\diff v^{n-1}}
	+\underbrace{\sum_{k=1}^{n}\hat{a}_{n-k}^{(n)}\diff v^k}\quad \text{for $n\ge2$},\\
	&\,\hspace{3.4cm}J_{B2}^n\hspace{4.5cm} J_{L1}^n\nonumber
\end{align}
where the discrete kernels $\hat{a}_{n-k}^{(n)}$ are defined by 
\begin{align}	
	\hat{a}_{0}^{(1)}:=&\,a_{0}^{(1)},\quad
	\hat{a}_{0}^{(n)}:=\frac{1-\alpha}{2-\alpha}a_{0}^{(n)}
	+\frac{\eta^{(n)}_{1}}{r_{n}(1+r_{n})}	
	\quad\text{for $n\ge2$,}	\label{eq: FBDF2 kernels hat a 0}\\
	\hat{a}_{n-k}^{(n)}:=&\,a_{n-k}^{(n)}-\frac{\eta^{(n)}_{n-k}}{1+r_{k+1}}
	+\frac{\eta^{(n)}_{n-k+1}}{r_{k}(1+r_{k})}\quad\text{for $2\le k\le n-1$ $(n\ge3)$,}
	\label{eq: FBDF2 kernels hat a n-k}\\
	\hat{a}_{n-1}^{(n)}:=&\,a_{n-1}^{(n)}-\frac{\eta^{(n)}_{n-1}}{1+r_{2}}
	\qquad\qquad\text{for $n\ge2$}.	\label{eq: FBDF2 kernels hat a n-1}
\end{align} 
The first part $J_{B2}^n$ represents a variable-step BDF2-type formula 
\cite{LiLiao:2022,LiaoJiWangZhang:2021ch,LiaoTangZhou:2020bdf2,LiaoZhang:2021MC}.
The second part $J_{L1}^n$ represents the nonuniform L1-type formula 
\cite{LiaoMcLeanZhang:2019,LiaoZhuWang:2021} of Caputo derivative \eqref{def:Caputo derivative}. 
As seen, the discrete kernel $a_{0}^{(n)}$ is split into two parts with different 
weights $\frac{1}{2-\alpha}$ and $\frac{1-\alpha}{2-\alpha}$.
The technical reason of the above local-nonlocal splitting is that 
the local term $J_{B2}^n\rightarrow D_2v^n$ as $\alpha\rightarrow1^-$, while the convolution kernels $\hat{a}_{n-k}^{(n)}$ ($n\ge2$) of the nonlocal term $J_{L1}^n$ will vanish, 
$\hat{a}_{n-k}^{(n)}\rightarrow0$ as $\alpha\rightarrow1^-$.
So we can build a discrete  energy that is asymptotically compatible 
with the discrete energy of the BDF2 scheme for the CH model, 
see Remark \ref{remark:FBDF2 energy compatibility}.
It significantly updates the previous energy laws in \cite{JiZhuLiao:2022TFCH,LiaoZhuWang:2021,LiaoTangZhou:2020Anenergy}, 
in which the discrete energies are always not asymptotically compatible 
in the fractional order limit $\alpha\rightarrow1^-$ although 
the associated energy laws are asymptotically compatible.

\subsection{DGS of local part}
We shall handle the local term $J_{B2}^n$ by applying 
the recent step-scaling technique \cite{LiaoJiWangZhang:2021ch} for the BDF2 formula.
Next lemma presents an upper bound of step-ratio $r_k$ to ensure 
the DGS of $J_{B2}^n$.

\begin{lemma}\label{lem: upper bound r*}
	For any fractional index $\alpha\in(0,1)$, it holds that
	\begin{align*}
		\mathrm{g}(x,y,\alpha):=\frac{2+2(1+\alpha)x-\alpha x^{2-\frac{\alpha}{2}}}{1+x}
		-\frac{\alpha y^{2-\frac{\alpha}{2}}}{1+y}>0\quad\text{for $0\le x,y< r^*$,}
	\end{align*}
	where $r^*=r^*(\alpha)$ be the unique root of the equation $1/\alpha+(1+1/\alpha)r^*- (r^*)^{2-\frac{\alpha}{2}}=0$.
%	\begin{align}\label{eq: upper bound r*}
%		1/\alpha+(1+1/\alpha)r^*- (r^*)^{2-\frac{\alpha}{2}}=0\quad\text{for $\alpha\in(0,1)$.}
%	\end{align}	
\end{lemma}
\begin{proof}At first, we show the well-poseness of $r^*$ for $\alpha\in(0,1)$.
	Consider an auxiliary function 
	\begin{align}\label{eq: auxiliary function g1hat}
		g_1(z,\alpha):=1/\alpha+(1+1/\alpha)z-z^{2-\frac{\alpha}{2}}\quad\text{for $z>0$}
	\end{align}
	with $\partial_zg_1(z,\alpha):=1+1/\alpha-\brat{2-\alpha/2}z^{1-\frac{\alpha}{2}}$.
	Obviously, the equation $\partial_zg_1(z,\alpha)=0$ has a unique root $z_1>1$
	since $\partial_zg_1(1,\alpha)=1/\alpha+\alpha/2-1>0$. 
	Thus for any $\alpha\in(0,1)$, $g_1(z,\alpha)$ 
	is increasing in $(0,z_1)$ and decreasing in $(z_1,+\infty)$ with respect to $z$.
	We have 
	$$g_1(z_1,\alpha)=1/\alpha+\frac{1-\alpha/2}{2-\alpha/2}(1+1/\alpha)z_1>0\quad\text{and}\quad g_1(+\infty,\alpha)<0.$$
	Thus there exists a unique $r^*\in(z_1,+\infty)$ such that $g_1(r^*,\alpha)=0$ for any $\alpha\in(0,1)$.

	Differentiating the function ${\mathrm{g}}$ with respect to $x$ gives $\partial_x{\mathrm{g}}=\frac{\alpha}{2(1+x)^2}g_2(x,\alpha)$	
	with an auxiliary function $g_2(x,\alpha):=4-\brat{4-\alpha}x^{1-\frac{\alpha}{2}}
	-\brat{2-\alpha}x^{2-\frac{\alpha}{2}}$.
%\begin{align*}		
%	\partial_x{\mathrm{g}}=\frac{\alpha g_2(x,\alpha)}{2(1+x)^2}
%	\quad\text{with}\quad g_2(x,\alpha):=4-\brat{4-\alpha}x^{1-\frac{\alpha}{2}}
%		-\brat{2-\alpha}x^{2-\frac{\alpha}{2}}.
%\end{align*}
Obviously, the function $g_2$ is decreasing with respect to $x>0$ and $g_2(0,\alpha)g_2(1,\alpha)<0$
such that $g_2(x,\alpha)$ has a unique zero point $x_1\in(0,1)$ for any $\alpha\in(0,1)$.
Thus, ${\mathrm{g}}$ is increasing in $(0, x_1)$ and decreasing in
$(x_1,+\infty)$ for $x>0$. Furthermore, ${\mathrm{g}}$ is decreasing with respect to $y$ due to 
the simple fact $\partial_y{\mathrm{g}}<0$. It holds that
\begin{align*}
	{\mathrm{g}}(x,y,\alpha)>&\,\min\{{\mathrm{g}}(0,r^*,\alpha),{\mathrm{g}}(r^*,r^*,\alpha)\}
	={\mathrm{g}}(r^*,r^*,\alpha)=\frac{2\alpha}{1+r^*}g_1(r^*,\alpha)=0,
\end{align*}
since $r^*>z_1>1>x_1$. This completes the proof.
\end{proof}

\begin{figure}[htb!]
		\centering
	\subfigure[$r^*=r^*(\alpha)$]{
		\includegraphics[width=2.1in]{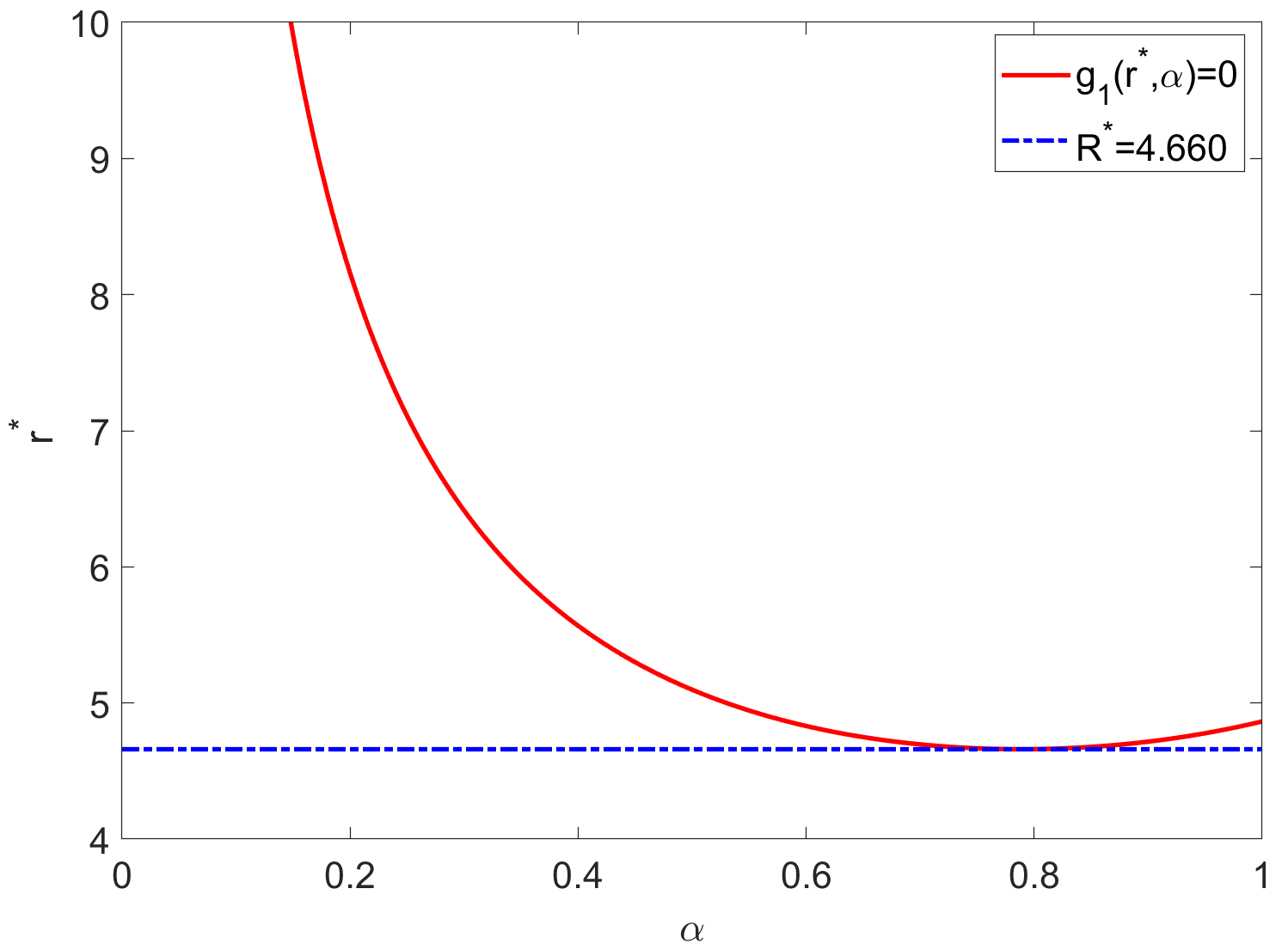}}
	\subfigure[$\gamma_{\max}=\gamma_{\max}(\alpha)$]{
		\includegraphics[width=2.1in]{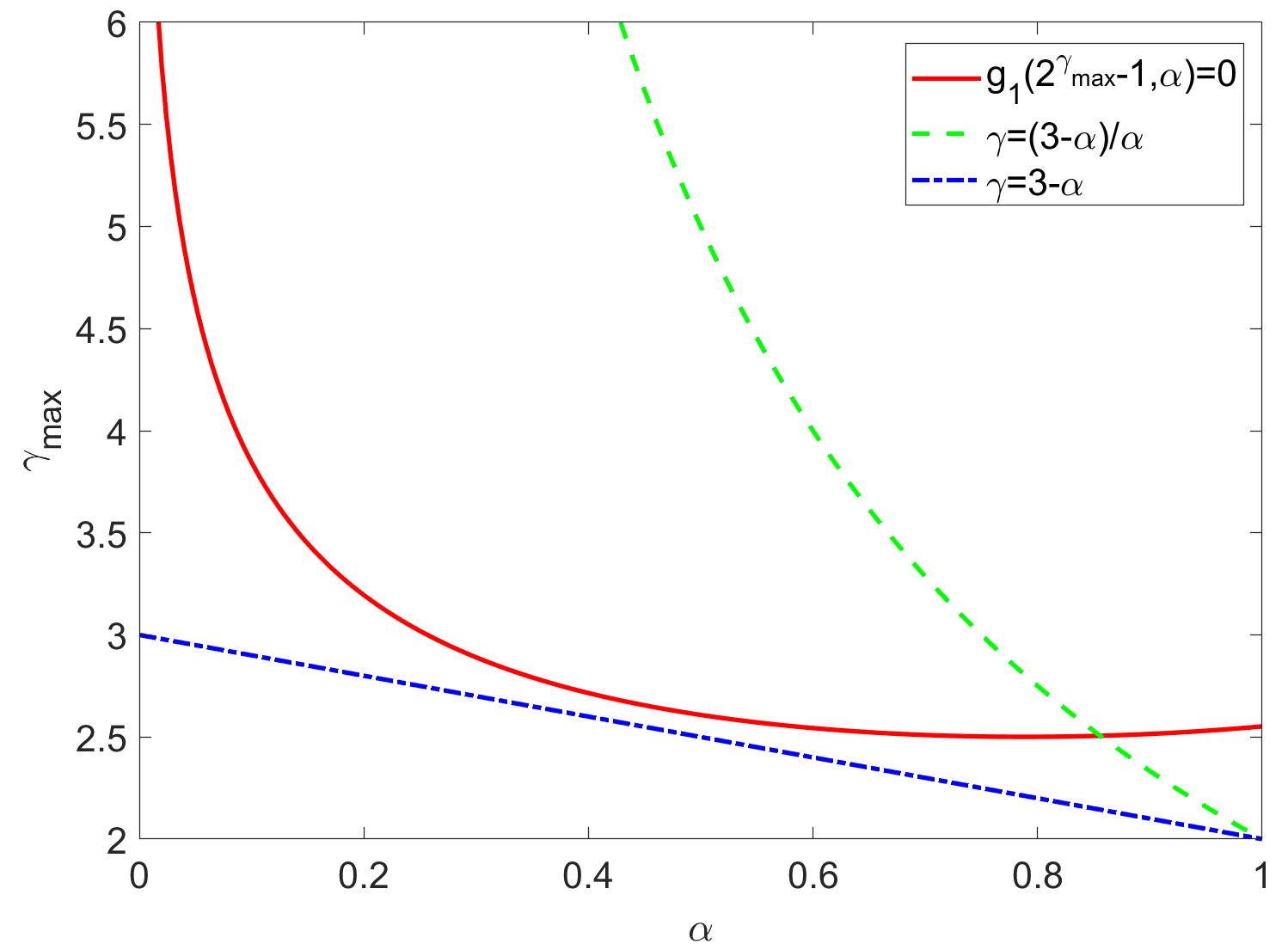}}
	\caption{Solution curves of $g_1(r^*,\alpha)=0$ and $g_1(2^{\gamma_{\max}}-1,\alpha)=0$.}
	\label{fig: implicitplotRemark1}
\end{figure}

\begin{remark}\label{remark: upper bound r*}
	It is easy to check that 
	$r^*(1)=r^\star\approx4.864$ and $r^*(\alpha)\to +\infty$ as $\alpha\to0^+$;
	but we are not able to find an explicit solution $r^*(\alpha)$ of $g_1(r^*,\alpha)=0$.		
	Differentiating this equation with respect to $\alpha$ gives	
	$$\frac{\zd r^*}{\zd\alpha}=\frac{\frac1{2}\alpha^2(r^*)^{2-\frac{\alpha}{2}}\ln r^*-1-r^*}{\brat{2-\frac{\alpha}{2}}\alpha^2(r^*)^{1-\frac{\alpha}{2}}-\alpha-\alpha^2}
	=\frac{(r^*)^{2-\frac{\alpha}{2}}g_3(r^*,\alpha)}
	{\brat{2-\frac{\alpha}{2}}\alpha(r^*)^{1-\frac{\alpha}{2}}-1-\alpha},$$
where the function $g_3(z,\alpha):=\frac1{2}\alpha\ln z-\frac{1+z}{1+(1+\alpha)z}$. 
Since the partial derivative $\partial_zg_3>0$ for $z>0$,
	it is not difficult to check that $g_3(z,\alpha)=0$ has a unique positive root for any  $\alpha\in(0,1)$. We solve the simultaneous equations of $g_1(r^*,\alpha)=0$ and $g_3(r^*,\alpha)=0$ numerically,
	and find a pair of approximate solutions  $(r^*,\alpha)\approx(4.660,0.7881)$. Let
	$R^{*}:=4.660.$ 
	It is checked that the function $g_1(R^*,\alpha)$ defined by \eqref{eq: auxiliary function g1hat} 
	has the minimum value,  that is, $g_1(R^*,\alpha)\ge g_1(R^*,0.7881)\approx6.67e-4$ for any $\alpha\in(0,1)$.
	%$$g_1(R^*,\alpha)\ge g_1(R^*,0.7881)\approx6.67e-4\quad\text{for any $\alpha\in(0,1)$.} $$
	That is to say, the solution of $g_1(r^*,\alpha)=0$ has a low bound $R^*$ 
	such that $r^*(\alpha)>R^*$ for any $\alpha\in(0,1)$, see Figure \ref{fig: implicitplotRemark1} (a),
	and the result of Lemma \ref{lem: upper bound r*} holds for $0\le x,y \le R^*=4.660$.	
\end{remark}

\begin{remark}\label{remark: graded meshes}
	In resolving the initial singularity of time-fractional problem \eqref{cont:TFCH model}, 
	the graded or quasi-graded meshes \cite{Kopteva:2019,Kopteva:2021,LiaoLiZhang:2018,LiaoMcLeanZhang:2019,
		LiaoMcLeanZhang:2021alikhanov,QuanWu:2022}
	would be practical to restore the convergence rate of our time-stepping scheme.  
	Due to the nonuniform setting, the FBDF2 method \eqref{def: FBDF2 formula} is 
	naturally applicable to these graded-like meshes; however, the step-ratio bound $r^*(\alpha)$ 
	in Lemma \ref{lem: upper bound r*} always restricts the value of graded parameter $\gamma$,
	which is adapted to the strength of the singularity. For 
	the standard graded mesh $t_k=T_0(k/N_0)^\gamma$ with a grading parameter $\gamma\ge1$, 
	the step-ratios 
	$r_k=\frac{k^\gamma-(k-1)^\gamma}{(k-1)^\gamma-(k-2)^\gamma}$
	is decreasing with respect to the level index $k$ with the maximum one $r_2=2^\gamma-1$.
	According to Remark \ref{remark: upper bound r*}, the maximum permissible parameter 
	$\gamma_{\max}=\gamma_{\max}(\alpha)$ can be determined by the equation $g_1(2^{\gamma_{\max}}-1,\alpha)=0$.
	Figure \ref{fig: implicitplotRemark1} (b) shows that $\gamma_{\max}(\alpha)>3-\alpha$ for any $\alpha\in(0,1)$.
	Thanks to  the theoretical results \cite[Theorem 4.1 and Remarks 4.2-4.4]{Kopteva:2021},
	it is adequate to attain the optimal convergence rate 
	at any positive time for the FBDF2 method \eqref{def: FBDF2 formula}.
	Nonetheless, the step-ratio bound $r^*(\alpha)$ in Lemma \ref{lem: upper bound r*} 
	is always inadequate to
	achieve the optimal rate of global convergence since 
	the maximum permissible parameter $\gamma_{\max}(\alpha)\not\ge (3-\alpha)/\alpha$
	on the standard graded mesh.
\end{remark}

\begin{lemma}\label{lem: DGS of JB2}
	Let $r_k< r^*(\alpha)$ for $k\ge2$. For the  local term $J_{B2}^n$ 
	in \eqref{eq: FBDF2-local nonlocal splitting}, it holds that
	\begin{align*}
		\brat{\diff v^{n}} J_{B2}^n\ge&\,
		\frac{\alpha r_{n+1}^{2-\frac{\alpha}{2}}}{2(1+r_{n+1})\tau_n^{\alpha}}
		\frac{ \brat{\diff v^{n}}^2}{\Gamma(3-\alpha)}
		-\frac{\alpha r_n^{2-\frac{\alpha}{2}}}{2(1+r_n)\tau_{n-1}^{\alpha}}
		\frac{\brat{\diff v^{n-1}}^2}{\Gamma(3-\alpha)}
		+\frac{\mathrm{g}(r_n,r_{n+1},\alpha)}{2\Gamma(3-\alpha)\tau_n^{\alpha}}\brat{\diff v^{n}}^2
	\end{align*}
for $n\ge2$, where the bound $r^*(\alpha)$ and $\mathrm{g}(x,y,\alpha)$ 
are defined by Lemma \ref{lem: upper bound r*}.
\end{lemma}
\begin{proof}
	%Recall that $\eta_0^{(n)}=\frac{\alpha\tau_n^{-\alpha}}{\Gamma(3-\alpha)}$. 
		We introduce the sequence $w_k$ via $\diff v^{k}:=\tau_k^{\alpha/2}w_k$
		%$$\diff v^{k}:=\tau_k^{\alpha/2}w_k\quad \text{for $k\ge1$} $$
		and derive from \eqref{eq: FBDF2-local nonlocal splitting} that
	\begin{align*}
		\brat{\diff v^{n}} J_{B2}^n=&\,\frac1{\Gamma(3-\alpha)}w_n^2+\frac{r_n}{1+r_{n}}\eta^{(n)}_{0}\tau_n^{\alpha}w_n^2
		-\frac{r_n^2}{1+r_{n}}\eta^{(n)}_{0}\tau_n^{\frac{\alpha}{2}}\tau_{n-1}^{\frac{\alpha}{2}}w_nw_{n-1}\\
		=&\,\frac{w_n^2}{\Gamma(3-\alpha)} +\frac{ \alpha r_n}{1+r_n}\frac{w_n^2}{\Gamma(3-\alpha)}
		-\frac{\alpha r_n^{2-\frac{\alpha}{2}}}{1+r_n}\frac{w_nw_{n-1}}{\Gamma(3-\alpha)}	\\
		\ge&\,\frac{1+(1+\alpha)r_n}{1+r_n}\frac{w_n^2}{\Gamma(3-\alpha)}
		-\frac{w_n^2}{\Gamma(3-\alpha)}\frac{\alpha r_n^{2-\frac{\alpha}{2}}}{2(1+r_n)}
		-\frac{w_{n-1}^2}{\Gamma(3-\alpha)}\frac{\alpha  r_n^{2-\frac{\alpha}{2}}}{2(1+r_n)}\\
		=&\,\frac{{\mathrm{g}}(r_n,r_{n+1},\alpha)}{2\Gamma(3-\alpha)}w_n^2
	+\frac{\alpha r_{n+1}^{2-\frac{\alpha}{2}}}{2(1+r_{n+1})}\frac{w_{n}^2}{\Gamma(3-\alpha)}
		-\frac{\alpha r_n^{2-\frac{\alpha}{2}}}{2(1+r_n)}\frac{w_{n-1}^2}{\Gamma(3-\alpha)}.
	\end{align*}
It leads to the claimed result by using $w_k=\tau_k^{-\alpha/2}\diff v^{k}$ for $k\ge1$.
\end{proof}

\subsection{DGS of nonlocal part}

To build up a DGS for the nonlocal term $J_{L1}^n$ 
in \eqref{eq: FBDF2-local nonlocal splitting}, we need
the following lemma, which is inspired by the proofs of \cite[Theorem 2.1]{JiZhuLiao:2022TFCH}
and \cite[Lemma 4.1]{QuanTangWangYang:2022}.

\begin{lemma}\label{lemma: minimum eigenvalue-Caputo kernels}
	For any fixed index $n\ge2$ and any positive kernels $\{\chi_{n-j}^{(n)}\}_{j=1}^n$, 
	define the following auxiliary kernels $\mathsf{a}_{0}^{(n)}:=(2-\sigma_{\min})\chi_{0}^{(n)}$
	and $\mathsf{a}_{n-j}^{(n)}:=\chi_{n-j}^{(n)}$ for $1\le j\le n-1$.
%\[\mathsf{a}_{0}^{(n)}:=(2-\sigma_{\min})\chi_{0}^{(n)}\quad\text{and}\quad \mathsf{a}_{n-j}^{(n)}:=\chi_{n-j}^{(n)}\quad\text{for $1\le j\le n-1$.}\]
	Assume that there exists a constant $\sigma_{\min}\in[0,2)$
	such that the modified kernels $\mathsf{a}_{n-j}^{(n)}$ satisfy
	\begin{description}	
		\item[(Row decrease)] \quad 
		$\displaystyle \mathsf{a}^{(n)}_{n-j-1}\ge \mathsf{a}^{(n)}_{n-j}>0$\; for $1\le j\le n-1$;	
		\item[(Column decrease)] \quad
		$\displaystyle \mathsf{a}^{(n-1)}_{n-1-j}\ge\mathsf{a}_{n-j}^{(n)}$\; for $1\le j\le n-1$;		
		\item[(Algebraic convexity)] \quad
		$\displaystyle \mathsf{a}^{(n-1)}_{n-2-j}-\mathsf{a}^{(n-1)}_{n-1-j}\ge
		\mathsf{a}^{(n)}_{n-1-j}-\mathsf{a}^{(n)}_{n-j}$\; for $1\le j\le n-2$.
	\end{description}
Then for any real sequence $\{w_k\}_{k=1}^n$, 
the following discrete gradient structure holds,
	\begin{align*}
		2w_n\sum_{j=1}^n\chi_{n-j}^{(n)}w_j=Y[\vec{w}_n]-Y[\vec{w}_{n-1}]+
		\sigma_{\min}\chi_{0}^{(n)}w_n^2+Y_R[\vec{w}_n]\qquad
	\text{for $n\ge1$,}
	\end{align*}
	where $\vec{w}_n=(w_n,w_{n-1},\cdots,w_1)^T$, and the 
	nonnegative functionals $Y$ and $Y_R$ are  defined by
	\begin{align*}
		Y[\vec{w}_n]:=&\,\sum_{j=1}^{n-1}\brab{\mathsf{a}^{(n)}_{n-j-1}-\mathsf{a}^{(n)}_{n-j}}
		\braB{\sum_{\ell=j+1}^nw_\ell}^2+\mathsf{a}^{(n)}_{n-1}\braB{\sum_{\ell=1}^nw_\ell}^2
		\quad	\text{for $n\ge1$,}\\
		Y_R[\vec{w}_n]:=&\,\sum_{j=1}^{n-2}\brab{
			\mathsf{a}^{(n-1)}_{n-2-j}-\mathsf{a}^{(n-1)}_{n-1-j}-\mathsf{a}^{(n)}_{n-j-1}+\mathsf{a}^{(n)}_{n-j}}
		\braB{\sum_{\ell=j+1}^{n-1}w_\ell}^2
		+\brab{\mathsf{a}^{(n-1)}_{n-2}-\mathsf{a}^{(n)}_{n-1}}\braB{\sum_{\ell=1}^{n-1}w_\ell}^2
	\end{align*}
for $n\ge2$.	Then the discrete convolution kernels 
$\chi_{n-k}^{(n)}$ are positive definite in the sense that 
	\begin{align*}
		2\sum_{k=1}^nw_k\sum_{j=1}^k\chi_{k-j}^{(k)}w_j\geq
		Y[\vec{w}_n]+\sigma_{\min}\sum_{k=1}^{n}\chi_{0}^{(k)}w_k^2\quad \text{for $n\ge1$}.
	\end{align*}
\end{lemma}
\begin{proof}
	Fix~$n$ and let $u_k:=\sum_{\ell=1}^kw_{\ell}$ for $0\le k\le n$ 
	such that $u_0=0$ and $ \diff u_k=w_k$ for $1\le k\le n$.
	%\[u_0=0\quad\text{and}\quad \diff u_k=w_k\quad\text{for $1\le k\le n$.}\]
	Recalling the following fact 
	$\sum_{j=1}^{n-1}\brab{\mathsf{a}^{(n)}_{n-j-1}-\mathsf{a}^{(n)}_{n-j}}
	=\mathsf{a}^{(n)}_{0}-\mathsf{a}^{(n)}_{n-1},$
	we derive that
		\begin{align*}
			2w_n\sum_{k=1}^n\mathsf{a}^{(n)}_{n-k}w_k
			=&\,2(\diff u_n)\sum_{k=1}^n\mathsf{a}^{(n)}_{n-k}\diff u_k
			=2(\diff u_n)\kbraB{\mathsf{a}^{(n)}_{0}u_n
				-\sum_{j=1}^{n-1}\brab{\mathsf{a}^{(n)}_{n-j-1}-\mathsf{a}^{(n)}_{n-j}}u_j}\\
			=&\,2(\diff u_n)\kbraB{
				\sum_{j=1}^{n-1}\brab{\mathsf{a}^{(n)}_{n-j-1}-\mathsf{a}^{(n)}_{n-j}}(u_n-u_j)
			+\mathsf{a}^{(n)}_{n-1}u_n}.
		\end{align*}
		By using the fact $2\diff u_n(u_n-u_j)=(u_n-u_j)^2-(u_{n-1}-u_j)^2+(\diff u_n)^2,$
%		\begin{align*}
%			2\diff u_n(u_n-u_j)
%		%	=&\,2(u_n-u_j)\kbra{(u_n-u_j)-(u_{n-1}-u_j)}\\
%		=&\,(u_n-u_j)^2-(u_{n-1}-u_j)^2+(u_n-u_{n-1})^2\quad\text{for $0\le j\le n-1$,}
%	\end{align*}
		it follows that
		\begin{align*}
			2w_n\sum_{k=1}^n\mathsf{a}^{(n)}_{n-k}w_k
			=&\,\sum_{j=1}^{n-1}\brab{\mathsf{a}^{(n)}_{n-j-1}-\mathsf{a}^{(n)}_{n-j}}(u_n-u_j)^2
				+\mathsf{a}^{(n)}_{n-1}u_n^2\\
				&\,-\sum_{j=1}^{n-1}\brab{\mathsf{a}^{(n)}_{n-j-1}-\mathsf{a}^{(n)}_{n-j}}(u_{n-1}-u_j)^2
				-\mathsf{a}^{(n)}_{n-1}u_{n-1}^2+\mathsf{a}^{(n)}_{0}(\diff u_n)^2\\
					=&\,\sum_{j=1}^{n-1}\brab{\mathsf{a}^{(n)}_{n-j-1}-\mathsf{a}^{(n)}_{n-j}}(u_n-u_j)^2
				+\mathsf{a}^{(n)}_{n-1}u_n^2\\
				&\,-\sum_{j=1}^{n-2}\brab{\mathsf{a}^{(n)}_{n-j-1}-\mathsf{a}^{(n)}_{n-j}}(u_{n-1}-u_j)^2
				-\mathsf{a}^{(n)}_{n-1}u_{n-1}^2+\mathsf{a}^{(n)}_{0}w_n^2.				
		\end{align*}
	By replacing the sequence $\{u_k\}$ with $\{w_k\}$, it is easy to find that
	\begin{align*}
		\mathsf{a}^{(n)}_{0}w_n^2+2w_n\sum_{k=1}^{n-1}\mathsf{a}^{(n)}_{n-k}w_k		
		=&\,Y[\vec{w}_n]-Y[\vec{w}_{n-1}]+Y_R[\vec{w}_n].
	\end{align*}
	The definition of $\mathsf{a}_{n-j}^{(n)}$ implies the claimed equality and
	completes the proof.
\end{proof}

Now we present the discrete gradient structure of the nonlocal part $J_{L1}^n$ 
in \eqref{eq: FBDF2-local nonlocal splitting}. 
\begin{theorem}\label{thm: DGS-nonlocal part}
	For the discrete  kernels $\hat{a}_{n-k}^{(n)}$ in 
	\eqref{eq: FBDF2 kernels hat a 0}-\eqref{eq: FBDF2 kernels hat a n-1}, 
	define the following auxiliary kernels 
	\begin{align}\label{eq: FBDF2 modified kernels A}
		A_{0}^{(n)}:=2\hat{a}_{0}^{(n)}\quad\text{and}\quad 
		A_{n-k}^{(n)}:=\hat{a}_{n-k}^{(n)}\quad\text{for $1\le k\le n-1$.}
	\end{align} 
	If the time-step ratios $r_k$ fulfill % $r_k\ge R_*$ for $k\ge2$.
\begin{align*}
	r_k\ge R_*=\frac{1}{9}\sqrt[3]{\frac{1}{2}\brat{189-81\sqrt{5}}}
	+\frac1{3}\sqrt[3]{\frac{1}{2}(7+3 \sqrt{5})}-\frac{1}{3}\approx 0.4753\quad\text{for $k\ge2$,}
\end{align*}	
then the auxiliary convolution kernels $A_{n-k}^{(n)}$ satisfy
	\begin{description}	
		\item[(a)] \;
		$ A^{(n)}_{n-k-1}> A^{(n)}_{n-k}>0$\; for $1\le k\le n-1$ $(n\ge2)$;	
		\item[(b)] \;
		$ A^{(n-1)}_{n-1-k}>A_{n-k}^{(n)}$\; for $1\le k\le n-1$ $(n\ge2)$;		
		\item[(c)]\;
		$ A^{(n-1)}_{n-2-k}-A^{(n-1)}_{n-1-k}
		> A^{(n)}_{n-k-1}-A^{(n)}_{n-k}$\; for $1\le k\le n-2$ $(n\ge3)$.
	\end{description}
	Then  for the  nonlocal term $J_{L1}^n$ 
	in \eqref{eq: FBDF2-local nonlocal splitting} with $n\ge2$, it holds that
	\begin{align*}
		(\diff v^n)\sum_{k=1}^n\hat{a}_{n-k}^{(n)}(\diff v^k)
		&\,\geq\frac12\kbra{\sum_{j=1}^{n-1}\brab{A^{(n)}_{n-j-1}-A^{(n)}_{n-j}}
		\braB{\sum_{\ell=j+1}^n\diff v^\ell}^2+A^{(n)}_{n-1}\braB{\sum_{\ell=1}^n\diff v^\ell}^2}\\
		&\,-\frac12\kbra{\sum_{j=1}^{n-2}\brab{A^{(n-1)}_{n-j-2}-A^{(n-1)}_{n-1-j}}
			\braB{\sum_{\ell=j+1}^{n-1}\diff v^\ell}^2+A^{(n-1)}_{n-2}\braB{\sum_{\ell=1}^{n-1}\diff v^\ell}^2}.
	\end{align*}
\end{theorem}
\begin{proof}
	Consider a function 
	$g_4(z):=3-\frac1{z^2(1+z)}$  with $g_4'(z)=\frac{3z+2}{z^3 (z+1)^2}>0$ for $z>0$.
It is easy to check that $R_*$ is the unique positive root of $g_4(z)=0$. If $r_k\ge R_*$ for $k\ge2$,
	it holds that
	\begin{align}\label{ieq: Ratio inequality nonlocal part}
		3-\frac{1}{r_{k}^2(1+r_{k})}\ge g_4(R_*)=0\quad\text{for $2\le k\le n$.}
	\end{align} 

According to Lemma \ref{lemma: minimum eigenvalue-Caputo kernels} with $\sigma_{\min}=0$, 
it remains to verify the three properties (a)-(c) of 
the kernels $A_{n-k}^{(n)}$. The remainder of this proof is presented separately in Section 4.	
\end{proof}

By applying Lemma \ref{lem: DGS of JB2} and Theorem \ref{thm: DGS-nonlocal part}, 
we obtain the DGS of the variable-step FBDF2 formula  \eqref{def: FBDF2 formula} 
via the splitting equality \eqref{eq: FBDF2-local nonlocal splitting} 
and the definitions \eqref{eq: FBDF2 kernels hat a 0}-\eqref{eq: FBDF2 kernels hat a n-1}.

\begin{theorem}\label{thm: DGS-FBDF2}
	For $n\ge1$, define a nonnegative functional $\mathcal{G}$ as follows,
	\begin{align*}
		\mathcal{G}[w_n]
		:=\frac{\alpha r_{n+1}^{2-\frac{\alpha}{2}}}{2(1+r_{n+1})\tau_n^{\alpha}}
		\frac{w_n^2}{\Gamma(3-\alpha)}
		+\frac12\sum_{j=1}^{n-1}\brab{A^{(n)}_{n-j-1}-A^{(n)}_{n-j}}
		\braB{\sum_{\ell=j+1}^nw_\ell}^2+\frac12A^{(n)}_{n-1}\braB{\sum_{\ell=1}^nw_\ell}^2.
	\end{align*}
	Assume that the time-step ratios $r_k$ fulfill 
	\begin{align}\label{eq: Ratio condition FBDF2}
		R_*\le r_k< r^*(\alpha) \quad\text{for $k\ge2$,}
	\end{align}
where the lower bound $R_*$  and the upper bound $r^*(\alpha)>R^*\approx4.660$
are given in Theorem \ref{thm: DGS-nonlocal part} and Lemma \ref{lem: upper bound r*} 
together with Remark \ref{remark: upper bound r*}, respectively. 
Then it holds that
	\begin{align*}
		(\diff v^n)(\dfd{\alpha}v)^{n}\geq \mathcal{G}\kbrab{\diff v^n}-\mathcal{G}\kbrab{\diff v^{n-1}}
		+\frac{{\mathrm{g}}(r_n,r_{n+1},\alpha)}{2\Gamma(3-\alpha)\tau_n^{\alpha}}\brat{\diff v^n}^2
	\quad \text{for $n\ge2$.}
	\end{align*}
\end{theorem}

\section{The FBDF2 scheme and discrete energy law}
\setcounter{equation}{0}

Now we apply the FBDF2 formula \eqref{def: FBDF2 formula}  to time integration of
TFCH flow \eqref{cont:TFCH model} and construct the following time-discrete scheme 
\begin{align}\label{scheme: implicit FBDF2 TFCH}
	\bra{\partial_\tau^\alpha \phi}^n=\kappa\Delta \mu^n\quad\text{with}\quad
	\mu^n=f(\phi^n)-\epsilon^2\Delta \phi^n \quad\text{for $n\ge 1$.}
\end{align}
Note that, the numerical scheme and the subsequent analysis can be extended in a straightforward way to the
fully discrete numerical schemes with some appropriate spatial discretization
preserving the discrete integration-by-parts formulas, such as the Fourier pseudo-spectral method \cite{ChengWangWise:2019,ChengWangWiseYue:2016}.

Obviously, the FBDF2
scheme \eqref{scheme: implicit FBDF2 TFCH} preserves the volume.
Actually, by taking the inner product of \eqref{scheme: implicit FBDF2 TFCH} with 1, one applies
the Green's formula to find that $\myinner{\brat{\partial_{\tau}^{\alpha}\phi}^n,1}
=\kappa\myinner{\Delta \mu^n,1}=0$.
%\begin{align*}
%	\myinner{\brat{\partial_{\tau}^{\alpha}\phi}^n,1}
%	=\kappa\myinner{\Delta \mu^n,1}=0.
%\end{align*}
Then the simple induction yields that $\myinner{\phi^n,1}=\myinner{\phi^{n-1},1}$ for $n\ge1$
with the help of the discrete formulas \eqref{def: FBDF2 formula} and \eqref{def: FBDF2 formula-start}. 
Then we obtain $\myinnerb{\phi^n,1}=\myinnerb{\phi^{0},1}$ for $n\ge1$.

%%%%%%%%%%%%%%%%%%%%%%%%%%%%%%%%%%%%%%%%%%%%%%%%%%%%%%%%%%%%%%%%%%%%%%%%%
\begin{theorem}\label{thm:FBDF2 convexity solvability TFCH}
	Under the time-step restriction
	\begin{align}\label{ineq: solvability condition FBDF2}
		\tau_n\le \sqrt[\alpha]{\frac{4\epsilon^2}{\kappa}\frac{2-\alpha+2r_n}{(1+r_n)\Gamma(3-\alpha)}}
		\quad\text{for $n\ge1$,}
	\end{align}
	the variable-step FBDF2 scheme \eqref{scheme: implicit FBDF2 TFCH} is uniquely solvable.
\end{theorem}
%%%%%%%%%%%%%%%%%%%%%%%%%%%%%%%%%%%%%%%%%%%%%%%%%%%%%%%%%%%%%%%%%%%%%%%%%
\begin{proof}
	For any fixed time-level index $n\ge1$,
	we consider the following energy functional $Q[z]$ on the space
	$\mathbb{V}^{*}:=\big\{z\in L^2\bra{\Omega}\,|\, \myinner{z,1}=\myinner{\phi^{n-1},1}\big\},$
	\begin{align*}
		Q[z]:=\frac{B_0^{(n)}}{2}\mynormb{z-\phi^{n-1}}_{-1}^2
		+\myinner{\mathcal{L}^{n-1},z-\phi^{n-1}}_{-1}
		+\frac{\kappa\epsilon^2}{2}\mynormb{\nabla z}^2
		+\frac{\kappa}4\mynormb{z}_{L^4}^4-\frac{\kappa}{2}\mynormb{z}^2,
	\end{align*}
	where we denote
	$\mathcal{L}^{n-1}:=\sum_{k=1}^{n-1}B_{n-k}^{(n)}\diff \phi^{k}$
	from the compact formulation \eqref{def: FBDF2 formula-compact}.
	With the help of the definitions \eqref{def: FBDF2 formula} and \eqref{def: FBDF2 formula-start}  
	together with \eqref{coeff: ank}-\eqref{coeff: eta nk}, 
	the restriction \eqref{ineq: solvability condition FBDF2}
	implies that 
	$$B_0^{(n)}\ge a_0^{(n)}+\frac{r_n}{1+r_n}\eta_0^{(n)}= \frac{2-\alpha+2r_n}{(1+r_n)\Gamma(3-\alpha)\tau_n^{\alpha}}\ge \frac{\kappa}{4\epsilon^2}\,.$$
	By using the generalized H\"{o}lder inequality
	$$\mynormb{v}^2\le \mynormb{\nabla v}\mynormb{v}_{-1}
	\le \epsilon^2\mynormb{\nabla v}^2
	+\frac{1}{4\epsilon^2}\mynormb{v}_{-1}^2\quad
	\text{for any $v\in\mathbb{\mathring V}$},$$
	we see that the energy functional $Q[z]$
	is convex with respect to $z$, that is,
	\begin{align*}
		\frac{\zd^2Q}{\zd s^2}[z+s\psi]\Big|_{s=0}
		&=B_0^{(n)}\mynormb{\psi}_{-1}^2
		+\kappa\epsilon^2\mynormb{\nabla \psi}^2
		-\kappa\mynormb{\psi}^2
		+3\kappa\mynormb{z\psi}^2\\
		&\ge \brab{B_0^{(n)}-\frac{\kappa}{4\epsilon^2}}\mynormb{\psi}_{-1}^2
		+3\kappa\mynormb{z\psi}^2>0.
	\end{align*}
	It is easy to show
	that the functional $Q[z]$ is coercive on $\mathbb{V}^*$, that is,
	\begin{align*}
		Q[z]&\ge \frac{\kappa}{4}\mynormb{z}_{L^4}^4-\frac{\kappa}{2}\mynormb{z}^2
		-\frac{1}{2B_0^{(n)}}\mynormb{\mathcal{L}^{n-1}}^2
		\ge\frac{\kappa}{2}\mynormb{z}^2
		-\frac{1}{2B_0^{(n)}}\mynormb{\mathcal{L}^{n-1}}^2
		-\kappa\abs{\Omega},
	\end{align*}
	where the inequality
	$\mynorm{v}_{L^4}^4\ge 4\mynorm{v}^2-4\abs{\Omega}$
	has been used in the last step.
	So the functional $Q[z]$ has a unique minimizer,
	which implies the FBDF2 scheme \eqref{scheme: implicit FBDF2 TFCH} exists a unique solution.
	%It completes the proof.
\end{proof}

\begin{remark}
	As $\alpha\rightarrow 1^-$, the variable-step
	FBDF2 scheme \eqref{scheme: implicit FBDF2 TFCH} approaches
	the BDF2 scheme
	\begin{align}\label{scheme: BDF2 CH}
		D_2\phi^n=\kappa\Delta \mu^n\quad\text{with}\quad
		\mu^n=\bra{\phi^n}^3-\phi^n-\epsilon^2\Delta \phi^n
		\quad\text{for $ n\ge 1$},
	\end{align}
	which is uniquely solvable if the time-step size
	$\tau_n\le \frac{4(1+2r_n)\epsilon^2}{(1+r_n)\kappa}$, 
	cf. \cite{LiaoJiWangZhang:2021ch,LiaoTangZhou:2020bdf2}.
	As seen, the time-step condition  \eqref{ineq: solvability condition FBDF2}
	is asymptotically compatible in the fractional order limit $\alpha\rightarrow 1^-$.
\end{remark}

Let $E\kbra{\phi^{n}}$ be the discrete version of the
free energy functional \eqref{cont:classical free energy},
\begin{align}\label{def:discrete original energy}
	E\kbra{\phi^{n}}:=\frac{\epsilon^2}{2}\mynormb{\nabla \phi^n}^2
	+\myinnerb{F(\phi^n),1}\quad\text{with}\quad
	F(\phi^n):=\frac14\brab{\brat{\phi^n}^2-1}^2\quad\text{for $n\ge 0$.}
\end{align}
With the nonnegative functional $\mathcal{G}$ in Theorem \ref{thm: DGS-FBDF2}, we define 
a modified energy functional
\begin{align}\label{def:discrete  energy}
	\mathcal{E}_\alpha\kbrab{\phi^{n}}:=E\kbrab{\phi^{n}}
	+\frac{1}{\kappa}\myinnerb{\mathcal{G}\kbrab{\diff \phi^n},1}_{-1}
	\quad\text{for $n\ge 1$.}
\end{align}

%%%%%%%%%%%%%%%%%%%%%%%%%%%%%%%%%%%%%%%%%%%%%%%%%%%%%%%%%%%%%%%%%%%%%%%%
\begin{theorem}\label{thm:FBDF2 energy law TFCH}
	Under the step-ratio constraint \eqref{eq: Ratio condition FBDF2} with the following step-size condition 
	\begin{align}\label{ineq: stable condition FBDF2}
		\tau_n\le \sqrt[\alpha]{\frac{4\epsilon^2}{\kappa}
			\frac{{\mathrm{g}}(r_n,r_{n+1},\alpha)}{\Gamma(3-\alpha)}}
		\quad\text{for $n\ge2$,}
	\end{align}
	the variable-step FBDF2 scheme \eqref{scheme: implicit FBDF2 TFCH}
	preserves the following discrete energy dissipation law
	\begin{align*}
		\partial_\tau\mathcal{E}_\alpha\kbra{\phi^n}\le 0 \quad \text{for $2\le n\le N$.}
	\end{align*}
\end{theorem}
%%%%%%%%%%%%%%%%%%%%%%%%%%%%%%%%%%%%%%%%%%%%%%%%%%%%%%%%%%%%%%%%%%%%%%%%

\begin{proof}
	Making the inner product of the equation \eqref{scheme: implicit FBDF2 TFCH} by  $\bra{-\Delta }^{-1}\diff \phi^{n}/\kappa$, one obtains
	\begin{align}\label{thmProof:FBDF2 energy law inner}
		\frac{1}{\kappa}\myinnerb{\bra{\partial_\tau^\alpha \phi}^n,\diff \phi^n}_{-1}
		+\myinnerb{\bra{\phi^n}^3-\phi^{n},\diff \phi^n}
		-\myinnerb{\epsilon^2\Delta \phi^n,\diff \phi^n}=0.
	\end{align}
	Recalling the simple inequality
	$4\brat{a^3-a}\bra{a-b}\ge \bra{a^2-1}^2-\bra{b^2-1}^2-2\bra{a-b}^2,$ 
	we can bound the second term of \eqref{thmProof:FBDF2 energy law inner} by
	\[
	\myinnerb{f\bra{\phi^n},\diff \phi^n}
	\ge \myinnerb{F(\phi^n),1}-\myinnerb{F(\phi^{n-1}),1}
	-\frac12\mynormb{\diff \phi^n}^2.
	\]
	Applying the identity $2a(a-b)=a^2-b^2+(a-b)^2$, we formulate the third term of 
	\eqref{thmProof:FBDF2 energy law inner} as 
	\begin{align*}
		-\myinnerb{\epsilon^2\Delta \phi^n,\diff \phi^n}
		=\frac{\epsilon^2}{2}\mynormb{\nabla \phi^n}^2
		-\frac{\epsilon^2}{2}\mynormb{\nabla \phi^{n-1}}^2
		+\frac{\epsilon^2}{2}\mynormb{\diff\nabla\phi^n}^2.
	\end{align*}
	Substituting the above results into \eqref{thmProof:FBDF2 energy law inner}, one has
	\begin{align}\label{thmProof:FBDF2 energy law norm}
		\frac{1}{\kappa}\myinnerb{\bra{\partial_\tau^\alpha \phi}^n,\diff \phi^n}_{-1}
		+\frac{\epsilon^2}{2}\mynormb{\diff \nabla \phi^n}^2
		-\frac12\mynormb{\diff \phi^n}^2
		+E\kbra{\phi^n}
		\le E\kbra{\phi^{n-1}}.
	\end{align}
	For the first term  of  \eqref{thmProof:FBDF2 energy law norm},
	Theorem \ref{thm: DGS-FBDF2} yields
	\begin{align*}
		\frac{1}{\kappa}\myinnerb{\bra{\partial_\tau^\alpha \phi}^n,\diff \phi^n}_{-1}
		\ge \frac{1}{\kappa}\myinnerb{\mathcal{G}\kbrab{\diff \phi^n},1}_{-1}-
		\frac{1}{\kappa}\myinnerb{\mathcal{G}\kbrab{\diff \phi^{n-1}},1}_{-1}
		+\frac{{\mathrm{g}}(r_n,r_{n+1},\alpha)}{2\kappa\Gamma(3-\alpha)\tau_n^{\alpha}}\mynormb{\diff \phi^n}_{-1}^2.
	\end{align*}
	By the generalized H\"{o}lder inequality, we have
	$$-\frac12\mynormb{\diff\phi^n}^2\ge -\frac{1}{8\epsilon^2}\mynormb{\diff\phi^n}_{-1}^2
	-\frac{\epsilon^2}{2}\mynormb{\diff\nabla \phi^n}^2.$$
	Inserting the above two estimates into \eqref{thmProof:FBDF2 energy law norm},
	one gets
	\begin{align*}
		\frac1{2\kappa}\braB{\frac{{\mathrm{g}}(r_n,r_{n+1},\alpha)}{\Gamma(3-\alpha)\tau_n^{\alpha}}
			-\frac{\kappa}{4\epsilon^2}}\mynormb{\diff\phi^n}_{-1}^2
		+\mathcal{E}_\alpha\kbrab{\phi^n}
		\le \mathcal{E}_\alpha\kbrab{\phi^{n-1}}.
	\end{align*}
	Then the claimed result follows from the time-step condition
	\eqref{ineq: stable condition FBDF2} immediately.
	%This completes the proof.
\end{proof}

\begin{remark}[asymptotical compatibility]\label{remark:FBDF2 energy compatibility}
	Under the time-step restriction, cf. \cite{LiaoJiWangZhang:2021ch},
	 $$\tau_n\le \frac{4\epsilon^2}{\kappa}\braB{\frac{2+4r_n-r_n^{3/2}}{1+r_n}
	 -\frac{r_{n+1}^{3/2}}{1+r_{n+1}}}
 =4{\mathrm{g}}(r_n,r_{n+1},1)\epsilon^2/\kappa\quad\text{for $n\ge2$,}$$
		the BDF2 scheme \eqref{scheme: BDF2 CH} for the classical CH model
	preserves the energy dissipation law,
	\begin{align}\label{inequality:BDF2 energy decay law}
		\partial_\tau\mathcal{E}\kbra{\phi^n}\le 0
		\quad \text{for  $2\le n\le N$,}
	\end{align}
where the modified energy
	\begin{align}\label{inequality:BDF2 energy}
	\mathcal{E}\kbra{\phi^n}:=E\kbra{\phi^n}
	+\frac{\sqrt{r_{n+1}}\tau_{n+1}}{2\kappa(1+r_{n+1})}\mynormb{\piff \phi^n}_{-1}^2
	\quad \text{for  $1\le n\le N$.}
\end{align}
	As the fractional index $\alpha\rightarrow 1^-$,
	the definition \eqref{eq: FBDF2 modified kernels A} together with 
	\eqref{eq: FBDF2 kernels hat a 0}-\eqref{eq: FBDF2 kernels hat a n-1} shows that
	the discrete kernels
	 $A_{n-k}^{(n)}\rightarrow 0$ for $1\le k\le n$. According to Theorem \ref{thm: DGS-FBDF2},
	 the discrete  energy \eqref{def:discrete  energy} is 
	 asymptotically compatible with the modified energy \eqref{inequality:BDF2 energy}, that is,
	 \[
	 \mathcal{E}_\alpha\kbra{\phi^n}\;\longrightarrow\;
	 \mathcal{E}\kbra{\phi^n}
	 \quad\text{as $\alpha\rightarrow 1^-$.}
	 \]
	Also, the energy dissipation law in Theorem \ref{thm:FBDF2 energy law TFCH}
	is asymptotically compatible with the energy dissipation law \eqref{inequality:BDF2 energy decay law}
	in the sense that
	\[
	\partial_\tau\mathcal{E}_\alpha\kbra{\phi^n}\le 0\;\longrightarrow\;
	\partial_\tau\mathcal{E}\kbra{\phi^n}\le 0
	\quad\text{as $\alpha\rightarrow 1^-$.}
	\]
	Additionally, it is not difficult to check that the modified energy \eqref{inequality:BDF2 energy}
	of BDF2 scheme and the discrete  energy \eqref{def:discrete  energy} 
	of the FBDF2 scheme degrade into the original Ginzburg-Landau 
	energy $E\kbra{\phi^n}$ in approaching the steady state ($\diff \phi^n\rightarrow0$), that is,
		\[
	\mathcal{E}_\alpha\kbra{\phi^n}\;\longrightarrow\;E\kbra{\phi^n}\quad\text{and}\quad
	\mathcal{E}\kbra{\phi^n}\;\longrightarrow\;E\kbra{\phi^n}
	\quad\text{as $t_n\rightarrow +\infty$.}
	\]
\end{remark}

\begin{remark}\label{remark:FBDF2 error analysis}
	To make the presentation more clear, the fully implicit time-stepping scheme 
	for the TFCH model	is only considered  in this paper; however, 
	Theorem \ref{thm: DGS-FBDF2} would be useful 
	to derive the discrete energy dissipation laws of some other 
	FBDF2 approximations, such as convex splitting scheme, stabilized scheme 
	and linearized schemes using the recent SAV techniques, 
	for other time-fractional phase-field models.
	Since the discrete convolution kernels $B_{k}^{(n)}$ in \eqref{def: FBDF2 formula-compact}
	lack the monotonicity \cite{Kopteva:2021,LvXu:2016} with respect to the subscript $k$, 
	another essential difficulty in theory is to 
	establish optimal $L^2$ norm error estimate of 
	the FBDF2 schemes on general nonuniform grids 
	with the step-ratio constraint \eqref{eq: Ratio condition FBDF2}.
	These issues will be studied and presented in separate reports. 
\end{remark}

\section{Proof of Theorem  \ref{thm: DGS-nonlocal part} (a)-(c)}
\setcounter{equation}{0}

%This section focuses on the theoretical properties  (a)-(c) of the modified convolution
%kernels~$A_{n-k}^{(n)}$ in \eqref{eq: FBDF2 modified kernels A} 
%on a general class of nonuniform meshes. 
From the definition \eqref{eq: FBDF2 modified kernels A} together 
with the discrete kernels $\hat{a}_{n-k}^{(n)}$ defined in 
\eqref{eq: FBDF2 kernels hat a 0}-\eqref{eq: FBDF2 kernels hat a n-1},  
we see that the discrete kernels $A_{n-k}^{(n)}$ are defined by  $a_{n-k}^{(n)}$ and 
$\eta^{(n)}_{n-k}$. Some technical lemmas on the positive coefficients $a_{n-k}^{(n)}$ and 
$\eta^{(n)}_{n-k}$ will be established firstly 
via their definitions \eqref{coeff: ank} and \eqref{coeff: eta nk}.

\begin{lemma}\label{lem: ank diff}
		For $n\ge2$, define the following ``bridging" integrals 
	$I^{(n)}_{n-k}$ and $J^{(n)}_{n-k}$ as follows,
	\begin{align}
		I_{n-k}^{(n)}:=&\,\int_{t_{k-1}}^{t_{k}}\frac{t-t_{k}}{\tau_{k}}\omega_{-\alpha}(t_n-t)\zd{t}
		\quad\text{for $1\leq k\leq n$,}\label{def: Ink notations}\\
		J_{n-k}^{(n)}:=&\,\int_{t_{k-1}}^{t_{k}}\frac{t_{k-1}-t}{\tau_{k}}
		\omega_{-\alpha}(t_n-t)\zd{ t}	\quad\text{for $1\leq k\leq n-1.$}\label{def: Jnk notations}
	\end{align}
	The positive coefficients $a^{(n)}_{n-k}$ in \eqref{coeff: ank} satisfy
	\begin{itemize}
		\item[(i)]
		$\displaystyle a_{n-k-1}^{(n)}-a_{n-k}^{(n)}
		=I_{n-k-1}^{(n)}+J_{n-k}^{(n)}$\; for $1\leq k\leq n-1$.	
		\item[(ii)] $\frac{2(1-\alpha)}{2-\alpha}a_{0}^{(n)}-a_{1}^{(n)}
		=\frac{\alpha}{2-\alpha}\omega_{1-\alpha}(\tau_n)+J_{1}^{(n)}$\,.
	\end{itemize}
\end{lemma}
\begin{proof} The original idea comes from \cite[Lemma 4.6]{LiaoMcLeanZhang:2021alikhanov}
	but we include the proof for completeness.
Applying the definition \eqref{coeff: ank},
	we exchange the order of integration to find
	\begin{align}\label{lemproof: ank diff1}
		a_{n-k-1}^{(n)}-&\,\omega_{1-\alpha}(t_n-t_k)
		=\int_{t_{k}}^{t_{k+1}}\frac{\omega_{1-\alpha}(t_n-s)-\omega_{1-\alpha}(t_n-t_k)}{\tau_{k+1}}\zd s\nonumber\\
		=&\,-\frac{1}{\tau_{k+1}}\int_{t_{k}}^{t_{k+1}}\int_{t_{k}}^{s}\omega_{-\alpha}(t_n-t)\zd{t}\zd s
		=\int_{t_{k}}^{t_{k+1}}\frac{t-t_{k+1}}{\tau_{k+1}}\omega_{-\alpha}(t_n-t)\zd{t}
		=I_{n-k-1}^{(n)}
	\end{align}
	for $0\leq k\leq n-1$, and similarly,
	\begin{align}\label{lemproof: ank diff2}
		\omega_{1-\alpha}(t_n-t_k)-a_{n-k}^{(n)}
		=&\,\int_{t_{k-1}}^{t_{k}}\frac{\omega_{1-\alpha}(t_n-t_k)-\omega_{1-\alpha}(t_n-s)}{\tau_{k}}\zd s
		=J_{n-k}^{(n)}
	\end{align}
	for $1\leq k\leq n-1$. Hence the claimed equality (i) is obtained by summing up the above two equalities \eqref{lemproof: ank diff1} and \eqref{lemproof: ank diff2} for $1\leq k\leq n-1$.
	Now taking the index $k=n-1$ in \eqref{lemproof: ank diff2} yields
	$\omega_{1-\alpha}(\tau_n)=a_{1}^{(n)}+J_{1}^{(n)}$. By using the simple fact  $a^{(n)}_0=\omega_{1-\alpha}(\tau_n)/\brat{1-\alpha}$, we have
	\begin{align*}
		\frac{2(1-\alpha)}{2-\alpha}a_{0}^{(n)}-a_{1}^{(n)}=&\,
		\frac{\alpha}{2-\alpha}a_{1}^{(n)}+\frac{2}{2-\alpha}J_{1}^{(n)}
		=\frac{\alpha}{2-\alpha}\omega_{1-\alpha}(\tau_n)+J_{1}^{(n)}.
	\end{align*}
	This confirms the result (ii) and completes the proof.
\end{proof}

\begin{lemma}\label{lem: decreasing property ank bnk}
	The positive coefficients $\eta^{(n)}_{n-k}$ in \eqref{coeff: eta nk} satisfy 
	\begin{itemize}
		\item[(i)]
		$  \eta^{(n)}_{n-k}	<\eta^{(n-1)}_{n-k-1}$\; for $1\leq k\leq n-1$.
		\item[(ii)]
		$ \displaystyle
		\eta^{(n)}_{n-k-1}>r_{k+1}\eta^{(n)}_{n-k}$\; for $1\leq k\leq n-1$ and		
		$$\eta^{(n)}_{n-k-1}-r_{k+1}\eta^{(n)}_{n-k}
		<\eta^{(n-1)}_{n-k-2}-r_{k+1}\eta^{(n-1)}_{n-k-1}\quad\text{for $1\leq k\leq n-2$}.$$
	\end{itemize}
\end{lemma}
\begin{proof}
	The integration by parts yields, also see \cite[Lemma 2.1]{LiaoMcLeanZhang:2021alikhanov},
	\begin{align}\label{coeff: eta nk 2}
		\eta_{n-k}^{(n)}=
		-\frac{1}{\tau_k^2}\int_{t_{k-1}}^{t_k}\brat{s-t_{k-1}}\brat{t_{k}-s}\omega_{-\alpha}(t_n-s)\zd{s}>0
		\quad \text{for $1\leq k\leq n$.}
	\end{align}
Then it is easy to check that
\begin{align*}
	\eta^{(n)}_{n-k}<-\frac{1}{\tau_{k}^2}\int_{t_{k-1}}^{t_{k}}
	\brat{s-t_{k-1}}\brat{t_{k}-s}\omega_{-\alpha}(t_{n-1}-s)\zd{s}=\eta^{(n-1)}_{n-k-1}
	\quad \text{for $1\le k\leq n-1$.}
\end{align*}
The result (i) is verified. 
Now consider an auxiliary function 
	\begin{align}\label{def: varpi funtion}
		\eta_{n,k}(z):=\frac{1}{\tau_{k}}\int_{t_{k-1}}^{t_{k-1}+z\tau_k}
		\frac{2s-2t_{k-1}-z\tau_k}{\tau_k}\omega_{1-\alpha}(t_{n}-s)\zd{s}\quad \text{for $1\le k\leq n$}
	\end{align}
	with the first and second derivatives
	\begin{align*}
		\eta_{n,k}'(z)=&\,z\omega_{1-\alpha}(t_n-t_{k-1}-z\tau_k)-
		\frac{1}{\tau_{k}}\int_{t_{k-1}}^{t_{k-1}+z\tau_k}
		\omega_{1-\alpha}(t_{n}-s)\zd{s}\quad \text{for $1\le k\leq n$,}\\
		\eta_{n,k}''(z)=&\,-z\tau_k\omega_{-\alpha}(t_n-t_{k-1}-z\tau_k)\quad \text{for $1\le k\leq n$.}	
	\end{align*}
	Notice that $\eta_{n,k}(1)=\eta^{(n)}_{n-k}$, $\eta_{n,k}(0)=0$ and $\eta_{n,k}'(0)=0$ for $1\le k\leq n$.
The Cauchy differential mean-value theorem says that there exist $z_{1k}$, $z_{2k}\in(0,1)$ such that
	\begin{align*}
		\frac{\eta^{(n)}_{n-k-1}}{\eta^{(n)}_{n-k}}
		&=\frac{\eta_{n,k+1}(1)}{\eta_{n,k}(1)}
		=\frac{\eta_{n,k+1}'(z_{1k})}{\eta_{n,k}'(z_{1k})}=\frac{\eta_{n,k+1}''(z_{2k})}{\eta_{n,k}''(z_{2k})}
		\\
		&=\frac{\tau_{k+1}}{\tau_k}\brab{\frac{t_n-t_{k-1}-z_{2k}\tau_k}
			{t_n-t_{k}-z_{2k}\tau_{k+1}}}^{1+\alpha}> r_{k+1}
		\quad \text{for $1\leq k\leq n-1$}.
	\end{align*}
	To process the proof, we introduce another auxiliary function
	\begin{align*}
		\varPi_{n,k}(z):=\eta_{n,k+1}(z)-r_{k+1}\eta_{n,k}(z)\quad \text{for $1\le k\leq n-1$}
	\end{align*}
such that $\varPi_{n,k}(1)=\eta^{(n)}_{n-k-1}-r_{k+1}\eta^{(n)}_{n-k}$.
Thanks to the Cauchy differential mean-value theorem, there exist $z_{3n}$, $z_{4n}\in(0,1)$ such that
\begin{align*}
\frac{\varPi_{n,k}(1)}{\varPi_{n-1,k}(1)}
	=&\,\frac{\varPi_{n,k}'(z_{3n})}{\varPi_{n-1,k}'(z_{3n})}
	=\frac{\varPi_{n,k}''(z_{4n})}{\varPi_{n-1,k}''(z_{4n})}
	\\
	&=\frac{-\omega_{-\alpha}(t_{n}-t_{k}-z_{4n}\tau_{k+1})+\omega_{-\alpha}(t_{n}-t_{k-1}-z_{4n}\tau_{k})}
	{-\omega_{-\alpha}(t_{n-1}-t_{k}-z_{4n}\tau_{k+1})+\omega_{-\alpha}(t_{n-1}-t_{k-1}-z_{4n}\tau_{k})}<1
\end{align*}
for $1\leq k\leq n-2$, where the last inequality comes from the convexity of function $t^{-1-\alpha}$. 
	Thus the two inequalities of (ii) are valid and the proof is completed.
\end{proof}

\begin{lemma}\label{lem: bnk Ink Jnk}
	For $n\ge2$, the positive coefficients $\eta^{(n)}_{n-k}$ in \eqref{coeff: eta nk} satisfy 
	\begin{itemize}
		\item[(i)]
		$\displaystyle I_{n-k}^{(n)}> \eta^{(n)}_{n-k}$\; for $1\leq k\leq n$, and
	$$I_{n-k}^{(n)}-\eta^{(n)}_{n-k}<I_{n-k-1}^{(n-1)}-\eta^{(n-1)}_{n-k-1}\quad\text{for $1\leq k\leq n-1$.}$$	
			\item[(ii)]
			$\displaystyle J_{n-k}^{(n)}> 3\eta^{(n)}_{n-k}$\; for $1\leq k\leq n-1$, and
			$$J_{n-k}^{(n)}-3\eta^{(n)}_{n-k}<J_{n-k-1}^{(n-1)}-3\eta^{(n-1)}_{n-k-1}
			\quad\text{for $1\leq k\leq n-2$.}$$	
	\end{itemize}
\end{lemma}
\begin{proof}
By using the definition \eqref{def: Ink notations} of $I_{n-k}^{(n)}$ and the formula 
\eqref{coeff: eta nk 2}, one gets
	\begin{align}\label{coeff: eta nk I nk}
		I_{n-k}^{(n)}-\eta_{n-k}^{(n)}=
		-\frac{1}{\tau_{k}^2}\int_{t_{k-1}}^{t_k}(t_{k}-t)^2\omega_{-\alpha}(t_n-t)\zd{t}>0
		\quad \text{for $1\leq k\leq n$},
	\end{align}
and
\begin{align*}
	I_{n-k}^{(n)}-\eta_{n-k}^{(n)}<
	-\frac{1}{\tau_{k}^2}\int_{t_{k-1}}^{t_k}(t_{k}-t)^2\omega_{-\alpha}(t_{n-1}-t)\zd{t}
	=I_{n-k-1}^{(n-1)}-\eta^{(n-1)}_{n-k-1}\quad \text{for $1\leq k\leq n-1$}.
\end{align*}
Thus the result (i) is verified. Moreover, the definition \eqref{def: Jnk notations} gives
	\begin{align*}
		J_{n-k}^{(n)}-3\eta_{n-k}^{(n)}=
		-\frac{3}{\tau_k^2}\int_{t_{k-1}}^{t_k}(t-t_{k-1})\brat{t-t_{k-1}-2\tau_k/3}\omega_{-\alpha}(t_n-t)\zd{t}>0
	\end{align*}
	for $1\leq k\leq n-1$ because $(t_n-t)^{-1-\alpha}$ is increasing with repsect to $t$
	and  $$\int_{t_{k-1}}^{t_k}(t-t_{k-1})\brat{t-t_{k-1}-2\tau_k/3}\zd{t}=\tau_k^3\int_0^1s(s-2/3)\zd{s}=0.$$ 
	In the similar way, 
	\begin{align*}
		J_{n-k}^{(n)}-3\eta_{n-k}^{(n)}<&\,
		-\frac{3}{\tau_k^2}\int_{t_{k-1}}^{t_k}(t-t_{k-1})\brat{t-t_{k-1}-2\tau_k/3}\omega_{-\alpha}(t_{n-1}-t)\zd{t}\\
		=&\,J_{n-k-1}^{(n-1)}-3\eta^{(n-1)}_{n-k-1}\quad \text{for $1\leq k\leq n-2$}.
	\end{align*}
They verify the result (ii). It completes the proof.
\end{proof}

\subsection{Proof of Theorem \ref{thm: DGS-nonlocal part} (a)}

It is in the position to show that the discrete kernels $A_{j}^{(n)}$ in 
\eqref{eq: FBDF2 modified kernels A} are positive and monotonically decreasing
with respect to $j$.

\begin{proof}[Proof of Theorem \ref{thm: DGS-nonlocal part} (a)]
	The definitions \eqref{eq: FBDF2 kernels hat a 0}-\eqref{eq: FBDF2 kernels hat a n-1} 
	of discrete convolution  kernels $\hat{a}_{n-k}^{(n)}$ will be used frequently. 
	With the definitions \eqref{coeff: ank}-\eqref{coeff: eta nk}, 
	one has $A_{0}^{(n)}>0$ and
	\begin{align}\label{eq: FBDF2 modified kernels A mainpart}
		A_{n-k}^{(n)}\ge&\,a_{n-k}^{(n)}-\frac1{1+r_{k+1}}\eta^{(n)}_{n-k}
		=\int_{t_{k-1}}^{t_k}\frac{t_{k+1}+t_k-2s}{(1+r_{k+1})\tau_k^2}\omega_{1-\alpha}(t_{n}-s)\zd{s}>0
	\end{align}
	for $1\le k\le n-1$ due to the simple fact $\int_{t_{k-1}}^{t_k}\brat{t_{k+1}+t_k-2s}\zd{s}>0.$	
	It remains to prove the decreasing property $A_{n-k-1}^{(n)}\ge A_{n-k}^{(n)}$ $(n\ge k+1\ge2)$
	by the four separate cases:
	(a.1) $k=1$ for $n=2$, (a.2) $k=1$ for $n\geq3$, (a.3) $2\leq k\leq n-2$ for $n\geq4$,
	and  (a.4) $k=n-1$ for	$n\geq3$.
	\begin{itemize}
		\item[\textbf{(a.1)}] \textbf{The case $k=1$  ($n=2$)}. Lemma \ref{lem: ank diff} (ii) gives
		\begin{align*}
			A_{0}^{(2)}-A_{1}^{(2)}=
			\frac{2-2\alpha}{2-\alpha}a^{(2)}_0-a_{1}^{(2)}+\frac{(2+r_2)\eta^{(2)}_{1}}{r_2(1+r_2)}
			=\frac{\alpha\omega_{1-\alpha}(\tau_2)}{2-\alpha}+J_{1}^{(2)}
			+\frac{(2+r_2)\eta^{(2)}_{1}}{r_2(1+r_2)}>0.
		\end{align*}		
		\item[\textbf{(a.2)}] \textbf{The case $k=1$ ($n\ge3$)}. 
	Applying Lemma \ref{lem: ank diff} (i), we reformulate the difference term
	$A_{n-2}^{(n)}-A_{n-1}^{(n)}$ into a linear combination with nonnegative coefficients, 
		\begin{align}\label{eq: modified kernels A difference n-2}
			A_{n-2}^{(n)}-A_{n-1}^{(n)}=&\,
			a_{n-2}^{(n)}-a_{n-1}^{(n)}+\frac{\eta_{n-1}^{(n)}}{r_2}-\frac{\eta_{n-2}^{(n)}}{1+r_3}
			=I_{n-2}^{(n)}+J_{n-1}^{(n)}+\frac{\eta_{n-1}^{(n)}}{r_2}-\frac{\eta_{n-2}^{(n)}}{1+r_3}\nonumber\\
			=&\,\brab{I_{n-2}^{(n)}-\eta_{n-2}^{(n)}}+\frac{r_3\eta_{n-2}^{(n)}}{1+r_3}
			+\brab{J_{n-1}^{(n)}-3\eta_{n-1}^{(n)}}+\frac{(1+3r_2)\eta_{n-1}^{(n)}}{r_2}.
		\end{align}
	Then Lemma \ref{lem: bnk Ink Jnk} (i)-(ii) yields $A_{n-2}^{(n)}>A_{n-1}^{(n)}$ directly.	
		\item[\textbf{(a.3)}] \textbf{The general cases $k=2,3,\cdots, n-2$  ($n\ge4$)}. 
		Lemma \ref{lem: ank diff} (i) gives
		\begin{align}
			A_{n-k-1}^{(n)}-A_{n-k}^{(n)}=&\,
			a_{n-k-1}^{(n)}-a_{n-k}^{(n)}+\frac{\eta^{(n)}_{n-k}}{r_{k+1}}
			-\frac{\eta_{n-k-1}^{(n)}}{1+r_{k+2}}
			-\frac{\eta^{(n)}_{n-k+1}}{r_{k}(1+r_{k})}\nonumber\\
%			=&\,I_{n-k-1}^{(n)}+J_{n-k}^{(n)}+\frac{\eta^{(n)}_{n-k}}{r_{k+1}}
%			-\frac{\eta_{n-k-1}^{(n)}}{1+r_{k+2}}
%			-\frac{\eta^{(n)}_{n-k+1}}{r_{k}(1+r_{k})}\nonumber\\
			=&\,\brab{I_{n-k-1}^{(n)}-\eta_{n-k-1}^{(n)}}+\frac{r_{k+2}}{1+r_{k+2}}\eta_{n-k-1}^{(n)}
			+\frac{\eta^{(n)}_{n-k}-r_k\eta^{(n)}_{n-k+1}}{r_{k}^2(1+r_{k})}\nonumber\\
			&\,+\brab{J_{n-k}^{(n)}-3\eta^{(n)}_{n-k}}
			+\braB{\frac{1}{r_{k+1}}+3-\frac{1}{r_{k}^2(1+r_{k})}}
			\eta^{(n)}_{n-k}.\label{eq: modified kernels A difference n-k-1}
		\end{align}	
	With Lemma \ref{lem: decreasing property ank bnk} (ii) and 
	Lemma \ref{lem: bnk Ink Jnk} (i)-(ii) together with the step-ratio 
		constraint \eqref{ieq: Ratio inequality nonlocal part}, the linear combination
		\eqref{eq: modified kernels A difference n-k-1} with nonnegative coefficients lead to 
		the desired inequality.		
	\item[\textbf{(a.4)}] \textbf{The case $k=n-1$  ($n\ge3$)}. 
	By using Lemma \ref{lem: ank diff} (ii), we have
	\begin{align}\label{eq: modified kernels A difference 0n}
		A_{0}^{(n)}-A_{1}^{(n)}=&\,\frac{2-2\alpha}{2-\alpha}a^{(n)}_0-a_{1}^{(n)}
		+\frac{(2+r_n)\eta^{(n)}_{1}}{r_n(1+r_n)}-\frac{\eta^{(n)}_{2}}{r_{n-1}(1+r_{n-1})}\nonumber\\
		=&\, \frac{\alpha\omega_{1-\alpha}(\tau_n)}{2-\alpha}
		+\brab{J_{1}^{(n)}-3\eta^{(n)}_{1}}
		+\frac{\eta^{(n)}_{1}-r_{n-1}\eta^{(n)}_{2}}{r_{n-1}^2(1+r_{n-1})}\nonumber\\
		&\,+\braB{\frac{2+r_n}{r_{n}(1+r_{n})} +3-\frac{1}{r_{n-1}^2(1+r_{n-1})}}\eta^{(n)}_{1}.
	\end{align}
	Thus Lemma \ref{lem: decreasing property ank bnk} (ii) and 
	Lemma \ref{lem: bnk Ink Jnk} (ii) together with the step-ratio 
	constraint \eqref{ieq: Ratio inequality nonlocal part} arrive at the desired inequality, 
	$A_{0}^{(n)}>A_{1}^{(n)}$.		
\end{itemize}
In summary, $A_{n-k-1}^{(n)}> A_{n-k}^{(n)}>0$ for $1\le k\le n-1$ such that
Theorem \ref{thm: DGS-nonlocal part} (a) is verified.
	\end{proof}

\subsection{Proof of Theorem \ref{thm: DGS-nonlocal part} (b)-(c)}

It is easy to check that
the inequalities of Theorem \ref{thm: DGS-nonlocal part} (b)-(c) are equivalent to 
\begin{align}\label{ieq: equivalent Theorem (b)-(c)}
	A^{(n-1)}_{0}-A^{(n)}_{1}
	> A^{(n-1)}_{1}-A^{(n)}_{2}>\cdots>
	A^{(n-1)}_{n-3}-A^{(n)}_{n-2}
	> A^{(n-1)}_{n-2}-A^{(n)}_{n-1}>0
	\quad\text{for $n\ge2$.}
\end{align}
Thus we shall prove Theorem  \ref{thm: DGS-nonlocal part} (b)-(c) by the following two technical steps:
\begin{itemize}[itemindent=0.8cm]
	\item[\textbf{Step 1}.] Check the last inequality,  
	$A^{(n)}_{n-1}<A^{(n-1)}_{n-2}$ for $n\ge2$.
	\item[\textbf{Step 2}.] 
	Verify the general cases of \eqref{ieq: equivalent Theorem (b)-(c)} via the equivalent inequalities,
	\begin{align}\label{ieq: equivalent Theorem (c)}
		A^{(n)}_{n-k-1}-A^{(n)}_{n-k} <A^{(n-1)}_{n-2-k}-A^{(n-1)}_{n-1-k}
		\quad\text{for $1\le k\le n-2$ $(n\ge3)$}.
	\end{align}
\end{itemize}

\begin{proof}[Proof of Theorem \ref{thm: DGS-nonlocal part} (b)-(c)]	
	By the formula \eqref{eq: FBDF2 modified kernels A mainpart}, one has	
	\begin{align*}
		A_{n-1}^{(n)}=a_{n-1}^{(n)}-\frac1{1+r_{2}}\eta^{(n)}_{n-1}
		<&\,\int_{t_{0}}^{t_1}\frac{t_{2}+t_1-2s}{(1+r_{2})\tau_1^2}
		\omega_{1-\alpha}(t_{n-1}-s)\zd{s}\nonumber\\
		=&\,a_{n-2}^{(n-1)}-\frac1{1+r_{2}}\eta^{(n-1)}_{n-2}\quad\text{for $n\ge2$.}
	\end{align*}
	According to the definitions \eqref{eq: FBDF2 modified kernels A} and
	\eqref{eq: FBDF2 kernels hat a n-1}, it leads to $A^{(n)}_{n-1}<A^{(n-1)}_{n-2}$ for $n\ge3$.	
	For $n=2$,	the definitions \eqref{eq: FBDF2 kernels hat a 0} 
	and \eqref{eq: FBDF2 kernels hat a n-1} yield that	$A_1^{(2)}<a_{0}^{(1)}-\frac1{1+r_{2}}\eta^{(1)}_{0}<a_0^{(1)}<A_0^{(1)}.$
	Thus we obtain that
	$$A^{(n)}_{n-1}<A^{(n-1)}_{n-2}\quad\text{for $n\ge2$.}$$
	
	To complete this proof, it remains to verify \eqref{ieq: equivalent Theorem (c)}.
	Recalling the three different expressions 
	\eqref{eq: modified kernels A difference n-2}-\eqref{eq: modified kernels A difference 0n} of 
	the difference term $A^{(n)}_{n-k-1}-A^{(n)}_{n-k}$, we will consider the three separate cases:
 (c.1) $k=1$ for $n\geq3$, (c.2) $2\leq k\leq n-3$ for $n\geq5$, and  (c.3) $k=n-2$ for	$n\geq4$.
 \begin{itemize}
  	\item[\textbf{(c.1)}] \textbf{The case $k=1$ ($n\ge3$)}. Applying the linear combination form
  	\eqref{eq: modified kernels A difference n-2},
  	Lemma \ref{lem: decreasing property ank bnk} (i)
  	and Lemma \ref{lem: bnk Ink Jnk} (i)-(ii), we have
 	\begin{align*}
 		A_{n-2}^{(n)}-A_{n-1}^{(n)}
 		<&\,\brab{I_{n-3}^{(n-1)}-\eta_{n-3}^{(n-1)}}+\frac{r_3\eta_{n-3}^{(n-1)}}{1+r_3} 		\\
 		&\,+\brab{J_{n-2}^{(n-1)}-3\eta_{n-2}^{(n-1)}}+\frac{(1+3r_2)\eta_{n-2}^{(n-1)}}{r_2}
 		=A_{n-3}^{(n-1)}-A_{n-2}^{(n-1)}.
 	\end{align*}
 	\item[\textbf{(c.2)}] \textbf{The general cases $k=2,3,\cdots, n-3$  ($n\ge5$)}. By using the equality 
 	\eqref{eq: modified kernels A difference n-k-1}, Lemma \ref{lem: decreasing property ank bnk} and 
 	Lemma \ref{lem: bnk Ink Jnk} (i)-(ii) together with the step-ratio 
 	constraint \eqref{ieq: Ratio inequality nonlocal part}, we derive that
 	\begin{align}\label{eq: modified kernels A column difference n-k-1}
 		A_{n-k-1}^{(n)}&\,-A_{n-k}^{(n)}<
 		\brab{I_{n-k-2}^{(n-1)}-\eta_{n-k-2}^{(n-1)}}
 		+\frac{r_{k+2}}{1+r_{k+2}}\eta_{n-k-2}^{(n-1)}\nonumber\\
 		&\,+\brab{J_{n-k-1}^{(n-1)}-3\eta^{(n-1)}_{n-k-1}}
 		+\frac{1}{r_{k}^2(1+r_{k})}\bra{\eta^{(n-1)}_{n-k-1}-r_k\eta^{(n-1)}_{n-k}}\nonumber\\	
 		&\,+\braB{\frac{1}{r_{k+1}}+3-\frac{1}{r_{k}^2(1+r_{k})}}
 		\eta^{(n-1)}_{n-k-1}\quad\text{for $2\le k\le n-2$.}
 	 	\end{align}	
  	According to the formula \eqref{eq: modified kernels A difference 0n}, 
  	the case $k=n-2$ in the right hand side of the above inequality does not equal to 
  	the difference term $A_{0}^{(n-1)}-A_{1}^{(n-1)}$ and will be handled in 
  	the next step \textbf{(c.3)}. At this moment, the equality 
  	\eqref{eq: modified kernels A difference n-k-1} implies that
  	$$A_{n-k-1}^{(n)}-A_{n-k}^{(n)}<A_{n-k-2}^{(n-1)}-A_{n-k-1}^{(n-1)}\quad\text{for $2\le k\le n-3$.}$$
 	\item[\textbf{(c.3)}] \textbf{The case $k=n-2$  ($n\ge4$)}. By considering the equality 
 	\eqref{eq: modified kernels A difference n-k-1} with $k=n-2$, 
 	we apply Lemma \ref{lem: decreasing property ank bnk}, 
 	Lemma \ref{lem: bnk Ink Jnk} (ii) and the 
 	condition \eqref{ieq: Ratio inequality nonlocal part} to the last three terms and obtain that
 	\begin{align*}
 		A_{1}^{(n)}-A_{2}^{(n)}<&\,
 		\brab{I_{1}^{(n)}-\eta_{1}^{(n)}}
 		+\frac{r_{n}}{1+r_{n}}\eta_{1}^{(n)}
 		+\brab{J_{1}^{(n-1)}-3\eta^{(n-1)}_{1}}\nonumber\\
 		&\,	+\frac{\eta^{(n-1)}_{1}-r_{n-2}\eta^{(n-1)}_{2}}{r_{n-2}^2(1+r_{n-2})}	
 		+\braB{\frac{1}{r_{n-1}}+3-\frac{1}{r_{n-2}^2(1+r_{n-2})}}
 		\eta^{(n-1)}_{1}\,.
 	\end{align*}
 	Replacing the index $n$ with $n-1$ in the formulation \eqref{eq: modified kernels A difference 0n}, one has
 	\begin{align*}
 		A_{0}^{(n-1)}-A_{1}^{(n-1)}=&\,\frac{\alpha\omega_{1-\alpha}(\tau_{n-1})}{2-\alpha}
 		+\brab{J_{1}^{(n-1)}-3\eta^{(n-1)}_{1}}
 		+\frac{\eta^{(n-1)}_{1}-r_{n-2}\eta^{(n-1)}_{2}}{r_{n-2}^2(1+r_{n-2})}\nonumber\\
 		&\, +\braB{\frac{1}{r_{n-1}(1+r_{n-1})}+\frac{1}{r_{n-1}} +3-\frac{1}{r_{n-2}^2(1+r_{n-2})}}\eta^{(n-1)}_{1}\,.
 	\end{align*}
 	Thus, by the equality \eqref{lemproof: ank diff1} with $k=n-2$, a simple subtraction yields 
 	\begin{align}\label{eq: modified kernels A column difference end}
 		A_{1}^{(n)}-A_{2}^{(n)}-&\,\brab{A_{0}^{(n-1)}-A_{1}^{(n-1)}}
 		<I_{1}^{(n)}-\frac{\eta_{1}^{(n)}}{1+r_{n}}
 		-\frac{\alpha\omega_{1-\alpha}(\tau_{n-1})}{2-\alpha}\nonumber\\
 		=&\,a_{1}^{(n)}-\omega_{1-\alpha}(\tau_n+\tau_{n-1})-\frac{\eta_{1}^{(n)}}{1+r_{n}}
 		-\frac{\alpha\omega_{1-\alpha}(\tau_{n-1})}{2-\alpha}\nonumber\\
 		=&-\frac{1}{\tau_{n-1}^{\alpha}\Gamma(2-\alpha)}\frac{g_5(r_n,\alpha)}{1+r_n},
 	\end{align} 
 where $g_5$ is an auxiliary function defined for $0<\alpha<1$ and $z\ge0$,
 \begin{align*}
 	g_5(z,\alpha):=&
 	\frac{\alpha}{2-\alpha}\kbra{(z+1)^{2-\alpha}-z^{2-\alpha}}
 	-\alpha(z+1)^{1-\alpha}+\frac{\alpha(1-\alpha)}{2-\alpha}(z+1)\,.
 \end{align*}
Note that, 	the first derivative of $g_5$ with respect to $z$ reads
\begin{align*}
	\partial_zg_5
	=&\alpha(1-\alpha)
	\int_0^{1}\kbra{(z+s)^{-\alpha}-(z+1)^{-\alpha}}\zd{s}
	+\frac{\alpha(1-\alpha)}{2-\alpha}>0\quad\text{for $z>0$}
\end{align*}
such that  $g_5(z,\alpha)>g_5(0,\alpha)=0$. Then the desired inequality follows from 
\eqref{eq: modified kernels A column difference end}.
 \end{itemize}
All of the inequalities in \eqref{ieq: equivalent Theorem (c)} are verified 
and the proof of Theorem \ref{thm: DGS-nonlocal part} (b)-(c) is completed.
\end{proof}

\section{Numerical experiments}
\setcounter{equation}{0}

Some numerical examples are presented to support the theoretical results.
Always, we use the Fourier pseudo-spectral method on the domain 
$\Omega=(0,2\pi)^2$ with $128\times128$ uniform grids. In solving the 
nonlinear system at each time level, a simple fixed-point iteration is adopted with the
tolerance of $10^{-12}$. The so-called SOE technique \cite{JZQZ-2017soe} is applied to speedup the FBDF2 formula \eqref{def: FBDF2 formula} with the absolute tolerance error $\varepsilon:=10^{-12}$ and the cut-off time $\Delta{t}:=\tau_1$.

\begin{figure}[htb!]
	\centering
	\subfigure[The fractional order $\alpha=0.4$]{
		\includegraphics[width=2.1in]{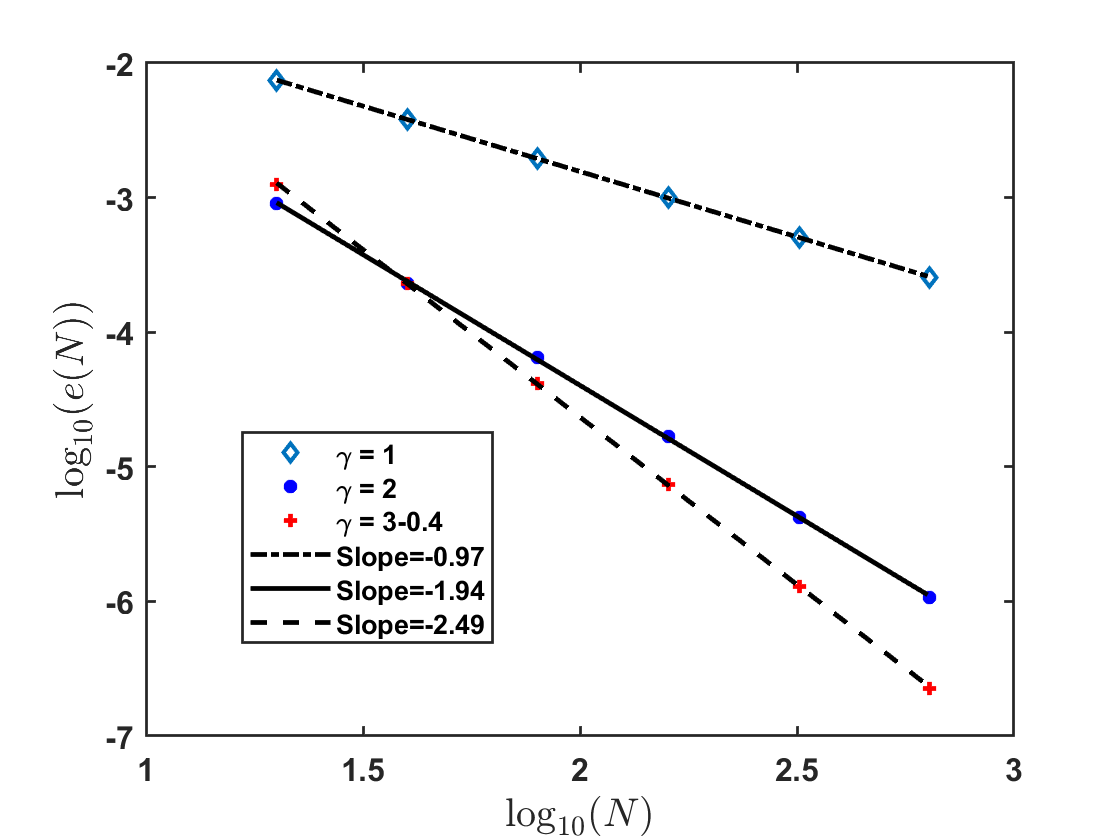}}
	\subfigure[The fractional order $\alpha=0.7$]{
		\includegraphics[width=2.1in]{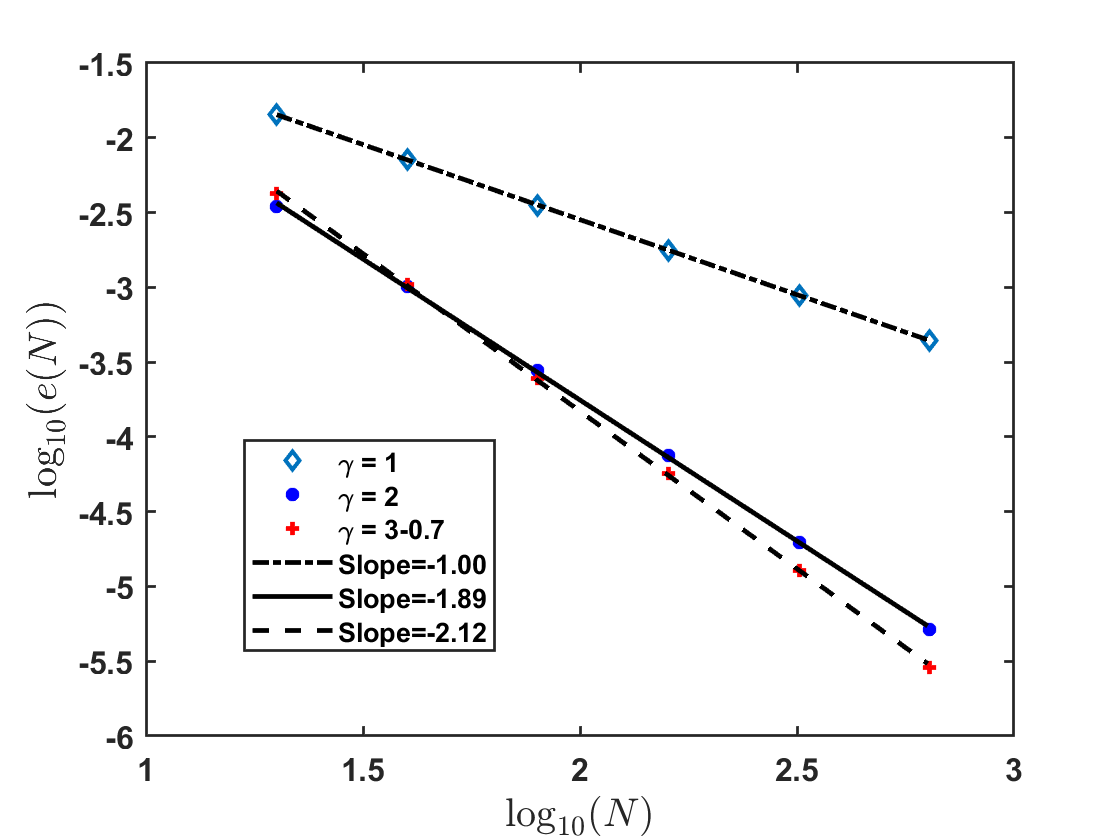}}
	\caption{Log-log plots of time accuracy for two different fractional orders. }
	\label{singu:L2 CH alpha 05}
\end{figure}

\begin{example}%\label{example:accuracy test}
	We consider the forced model 
	$\partial_{t}^{\alpha}\Phi=\kappa\Delta\mu-g$
	by adding a forcing term $g=g(\mathbf{x},t)$ to 
	the time-fractional Cahn-Hilliard model \eqref{cont:TFCH model} 
	on the domain $\mathbf{x}\in(0,2\pi)^{2}$. 
	Choose the model data $\kappa=1$, $\epsilon=0.5$ and $T=1$
	such that the exact solution $\Phi=\omega_{1+\alpha}(t)\sin x\sin y$.
\end{example}

We examine the temoral accuracy by considering
the graded mesh $t_k=(k/N)^\gamma$ for six different total number 
$N=20\times 2^{m-1}$ $(1\le m\le 6)$ such that
the temporal errors dominate the numerical error. 
In each run of the variable-step FBDF2 scheme \eqref{scheme: implicit FBDF2 TFCH}, 
the error $e(N):=\mynorm{\Phi^N-\phi^N}$ at the finial time $T=1$ is recorded. 
Figure \ref{singu:L2 CH alpha 05} depicts the log-log plots of the numerical errors 
for two different fractional orders $\alpha=0.4$ and $0.7$ with three grading parameters. 
As observed, the FBDF2 scheme \eqref{scheme: implicit FBDF2 TFCH} is accurate 
with the temporal order about $O(\tau^{\min\{\gamma,3-\alpha\}})$,
%(except the logarithmic factor $\ln N$ for $\gamma=3-\alpha$), 
which is in accordance with the previous results in \cite{Kopteva:2021,QuanWu:2022}.

\begin{example}\label{example:coarsen dynamic}
	Setting the mobility $\kappa=0.01$ and
	the interfacial thickness $\epsilon=0.05$,
	The coarsening dynamics of the time-fractional  Cahn-Hilliard model \eqref{cont:TFCH model}
	is simulated with a
	random initial field $\phi_0(\mathbf{x})=\textrm{rand}(\mathbf{x})$,
	where $\text{rand}(\mathbf{x})$ generates uniform random numbers
	between $-0.001$ to $0.001$.
\end{example}

\begin{figure}[htb!]
	\centering
	\subfigure[Original energy $E\kbra{\phi^n}$]{
		\includegraphics[width=2.1in]{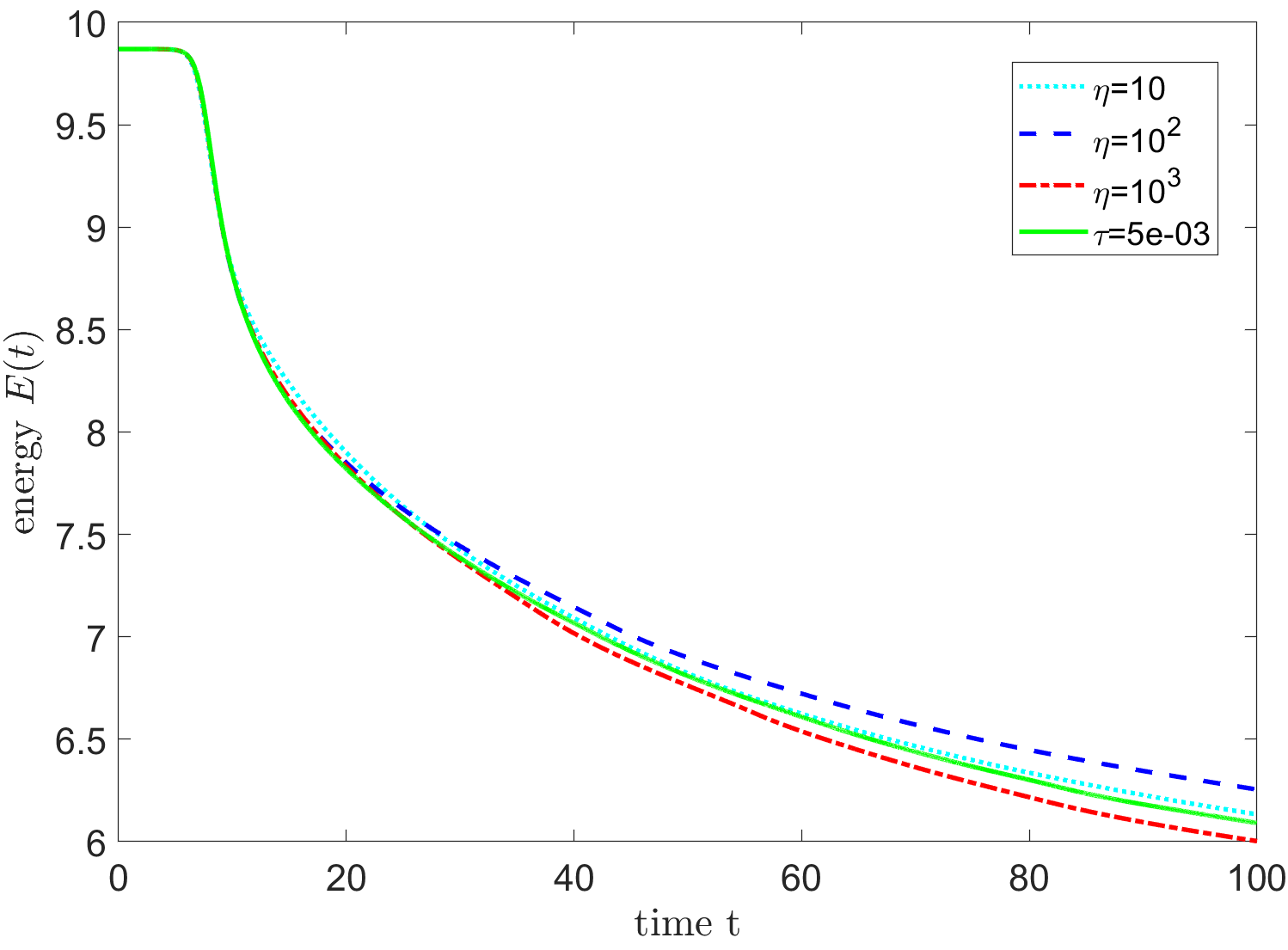}}
	\subfigure[Time steps $\tau_n$]{
		\includegraphics[width=2.1in]{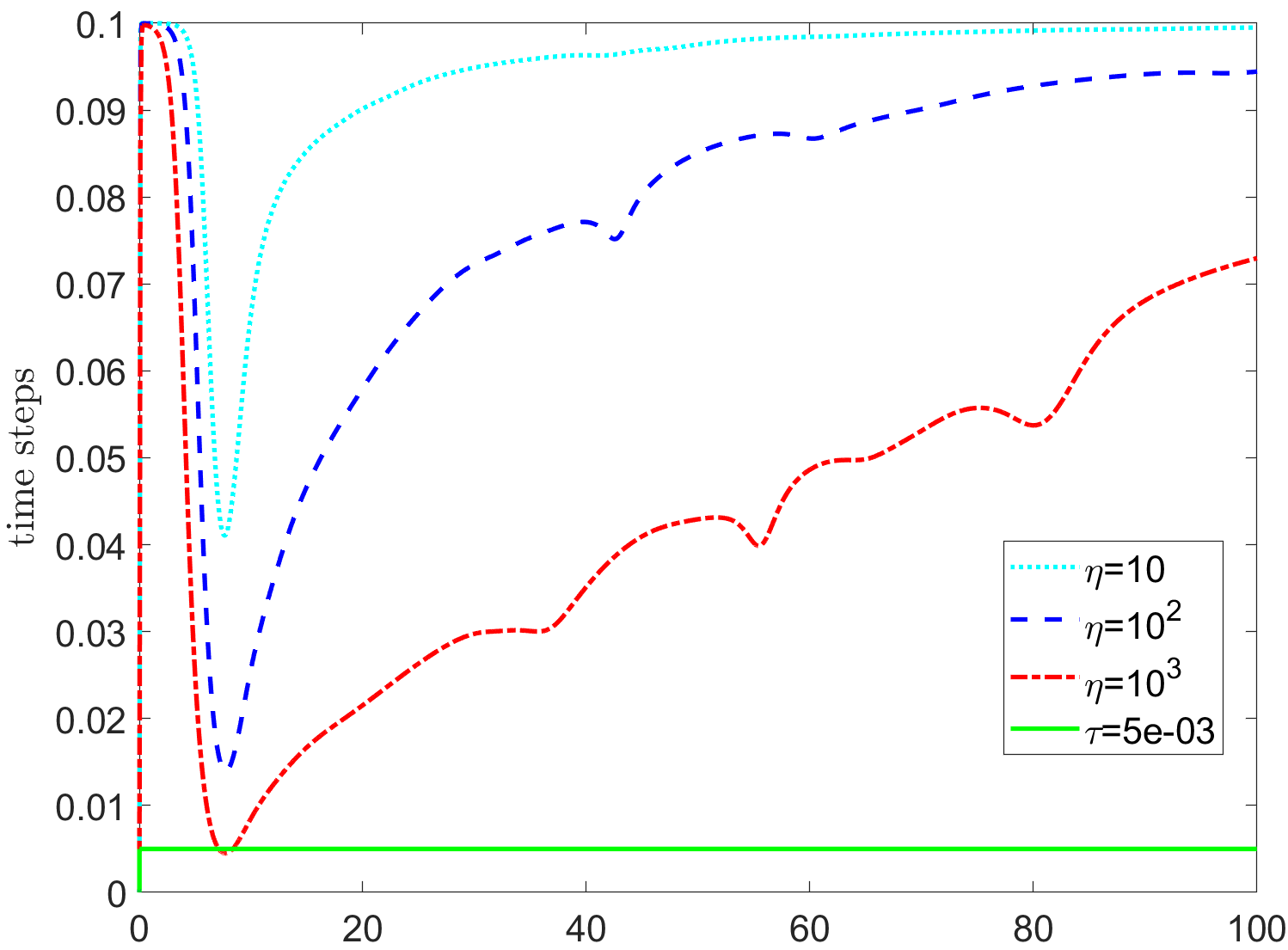}}
	\caption{Energy curves by uniform step and adaptive strategy with different parameters $\eta$.}
	\label{example:compar adap param}
\end{figure}

\begin{figure}[htb!]
	\centering
	\subfigure[Original energy $E\kbra{\phi^n}$]{
		\includegraphics[width=2.1in]{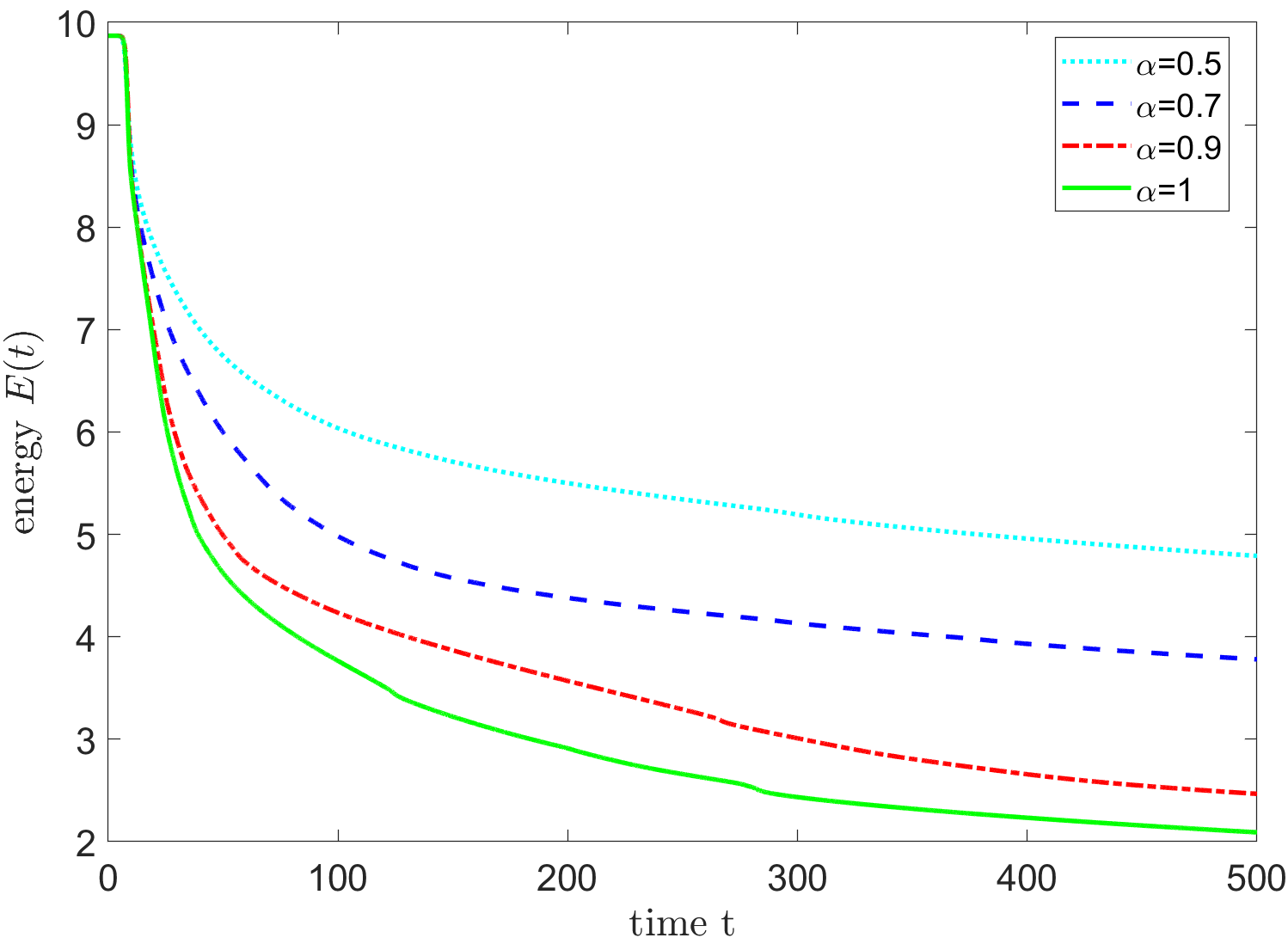}}
	\subfigure[Adaptive time steps $\tau_n$]{
		\includegraphics[width=2.1in]{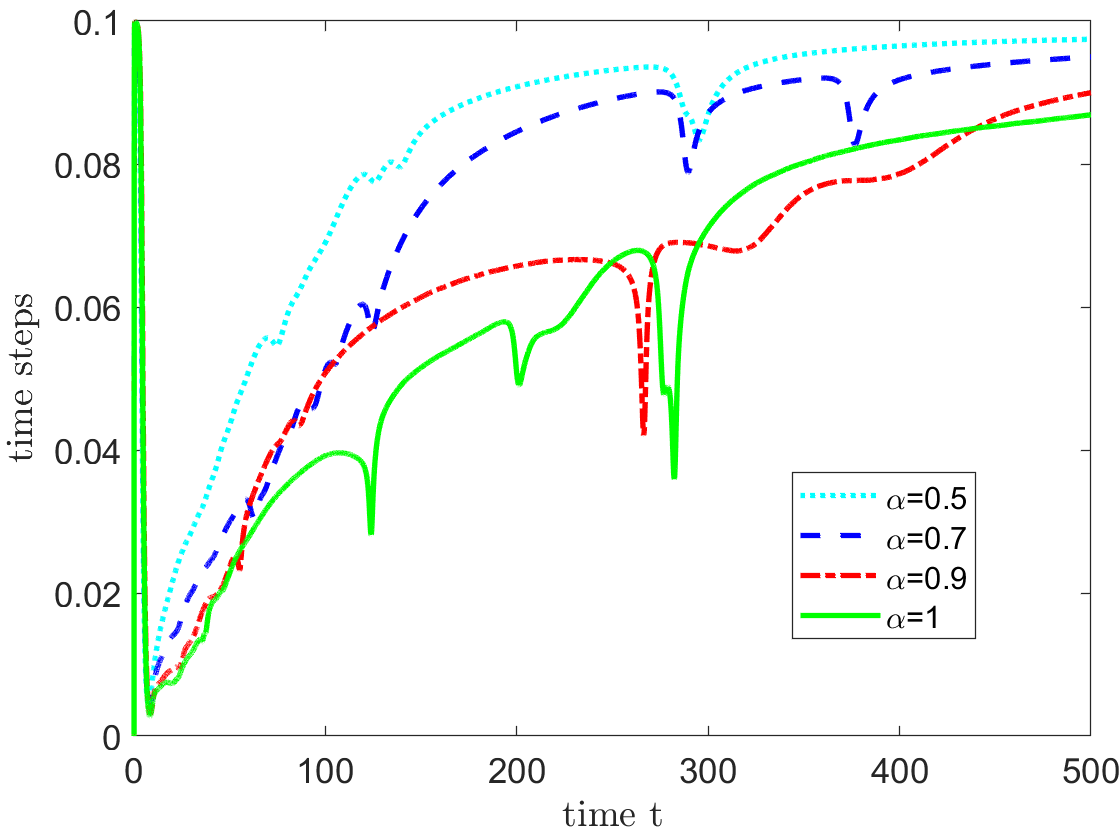}}
	\caption{Energy curves of  Example 2 with different fractional orders $\alpha$.}
	\label{example:FCH coarsen energy}
\end{figure}
\begin{figure}[htb!]
	\centering
	\subfigure[The fractional order $\alpha=0.2$, at time $t=10,30,50,100$]{
		\includegraphics[width=1.40in]{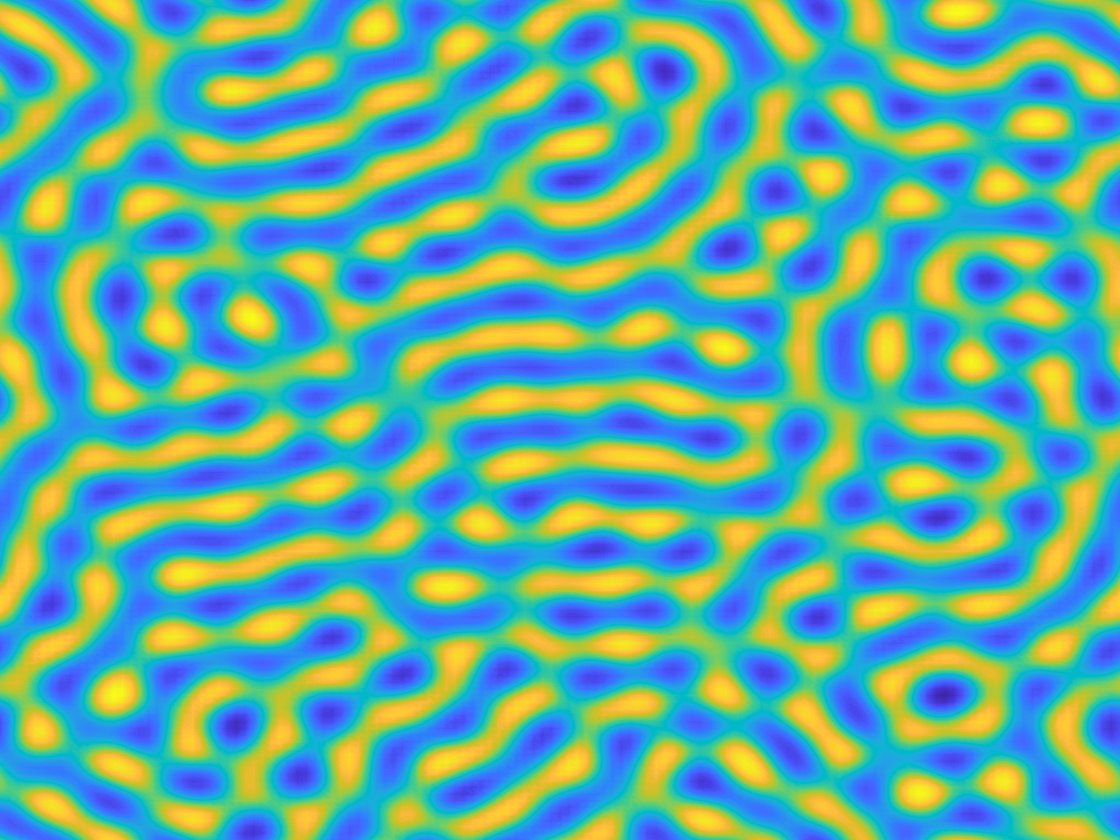}
		\includegraphics[width=1.40in]{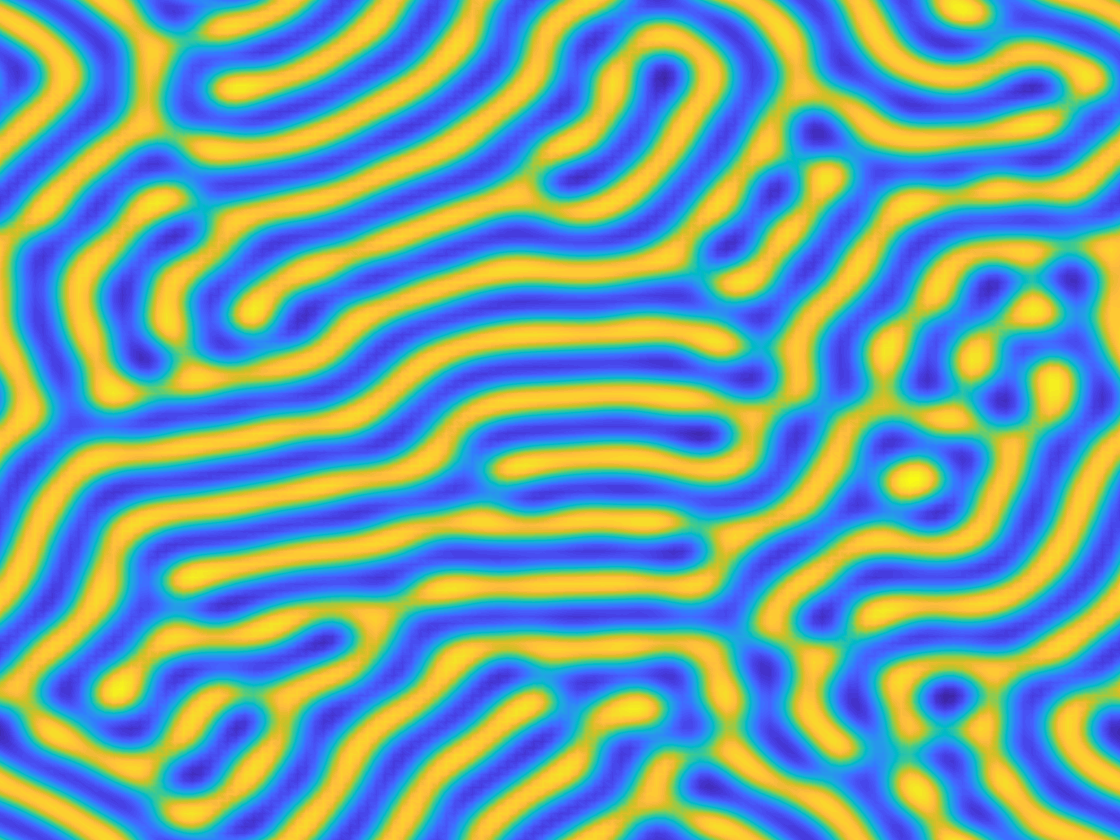}
		\includegraphics[width=1.40in]{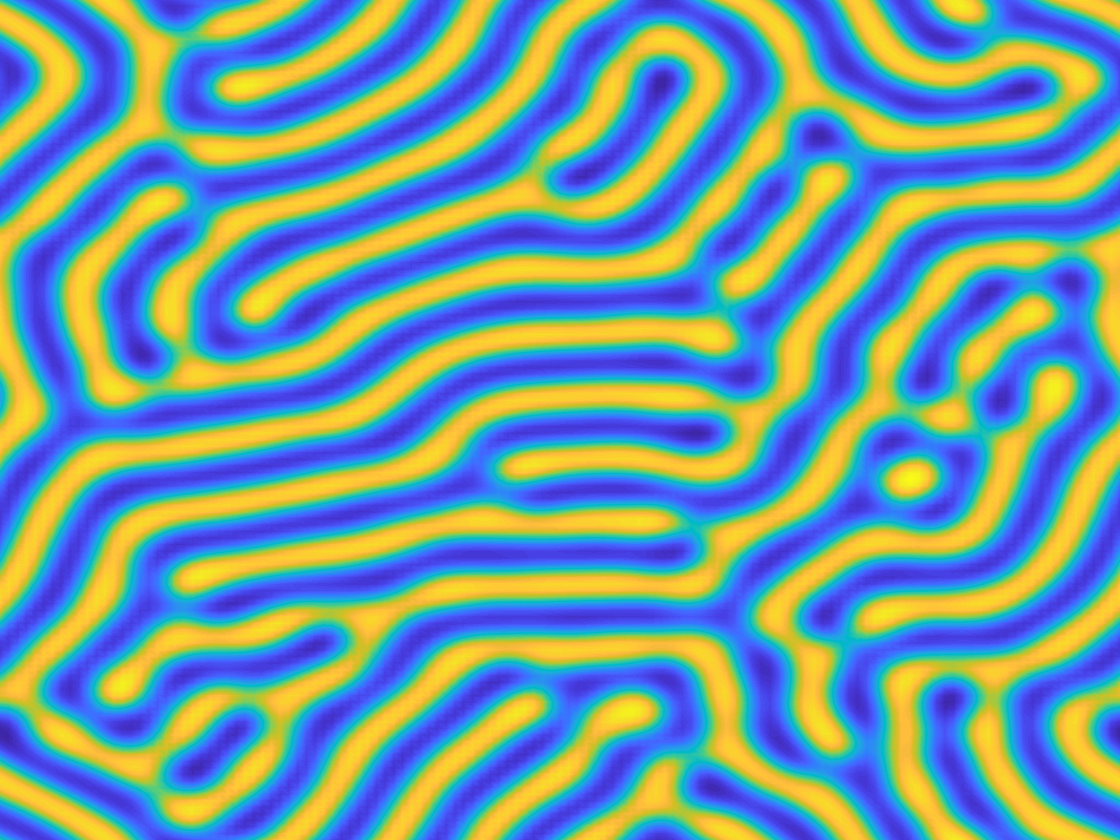}
		\includegraphics[width=1.40in]{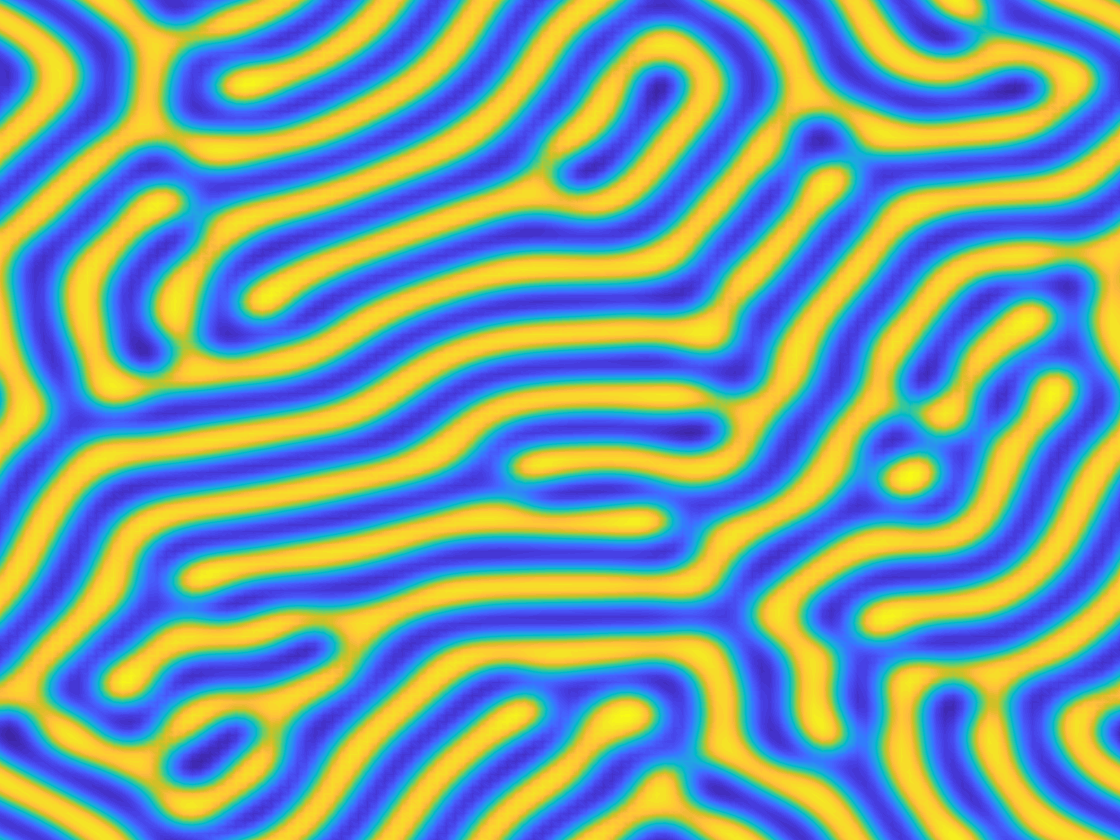}}
	\subfigure[The fractional order $\alpha=0.5$, at time $t=10,30,50,100$]{\includegraphics[width=1.40in]{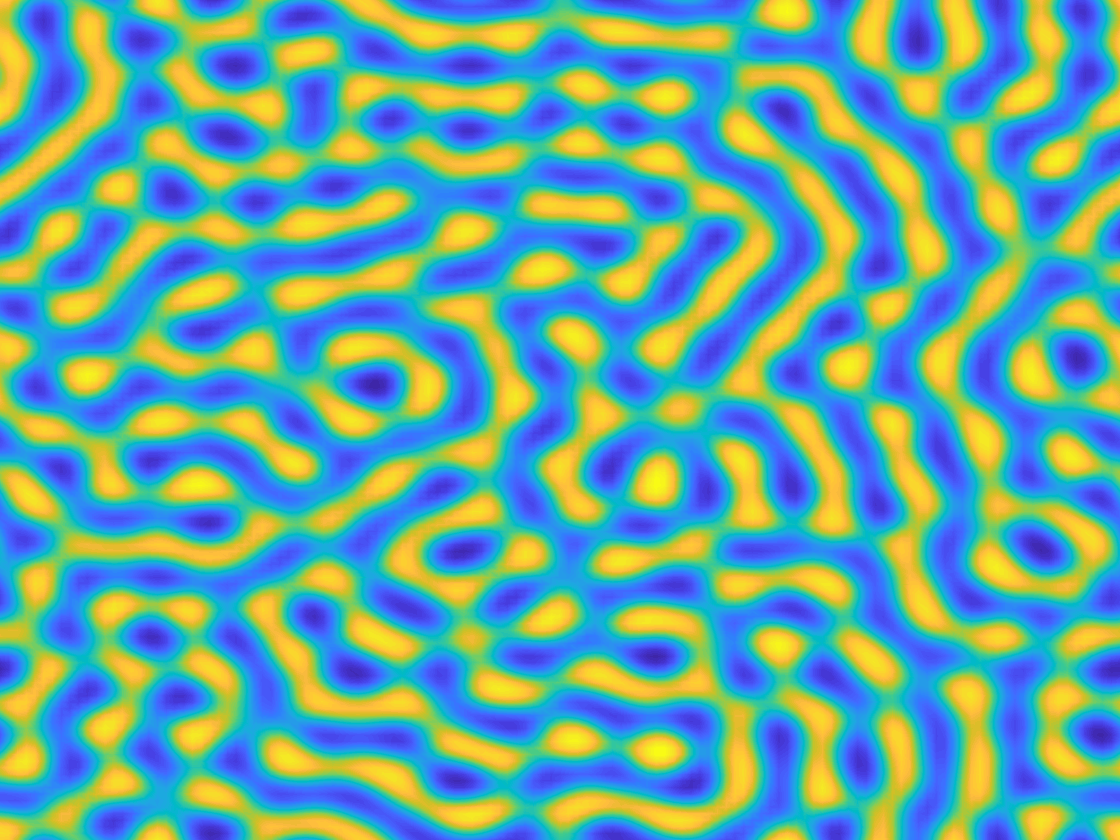}
		\includegraphics[width=1.40in]{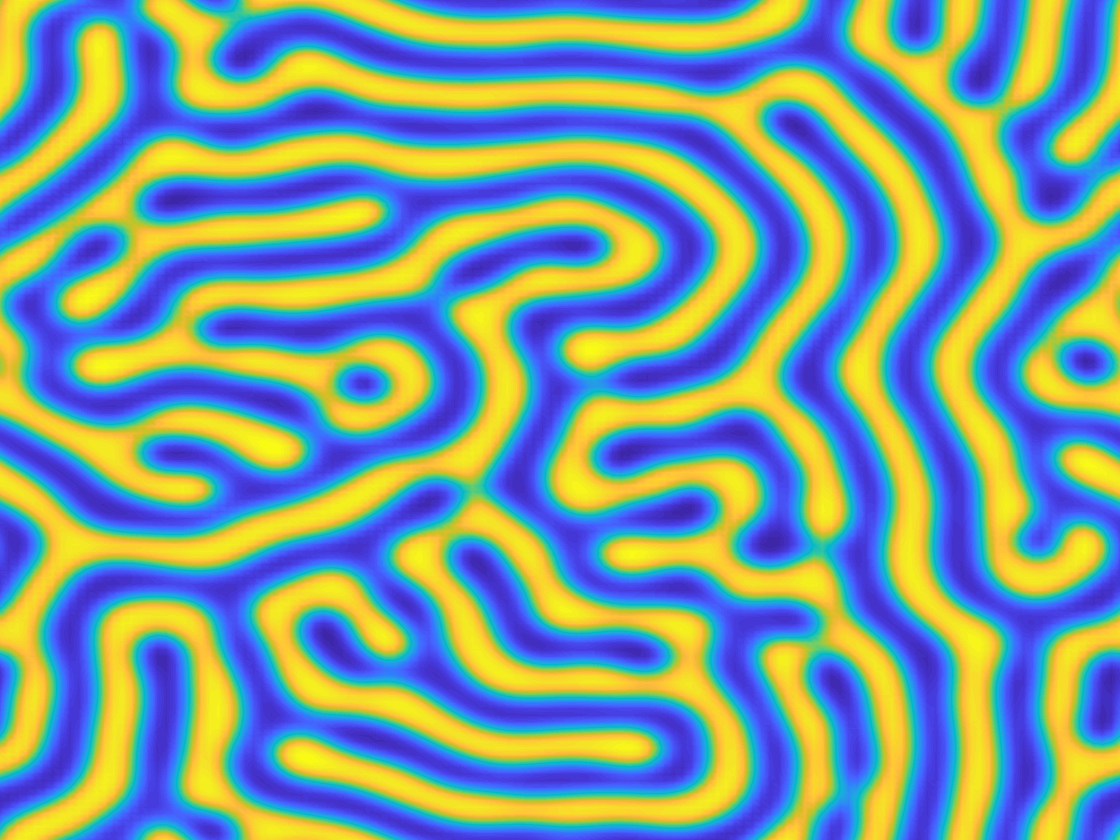}
		\includegraphics[width=1.40in]{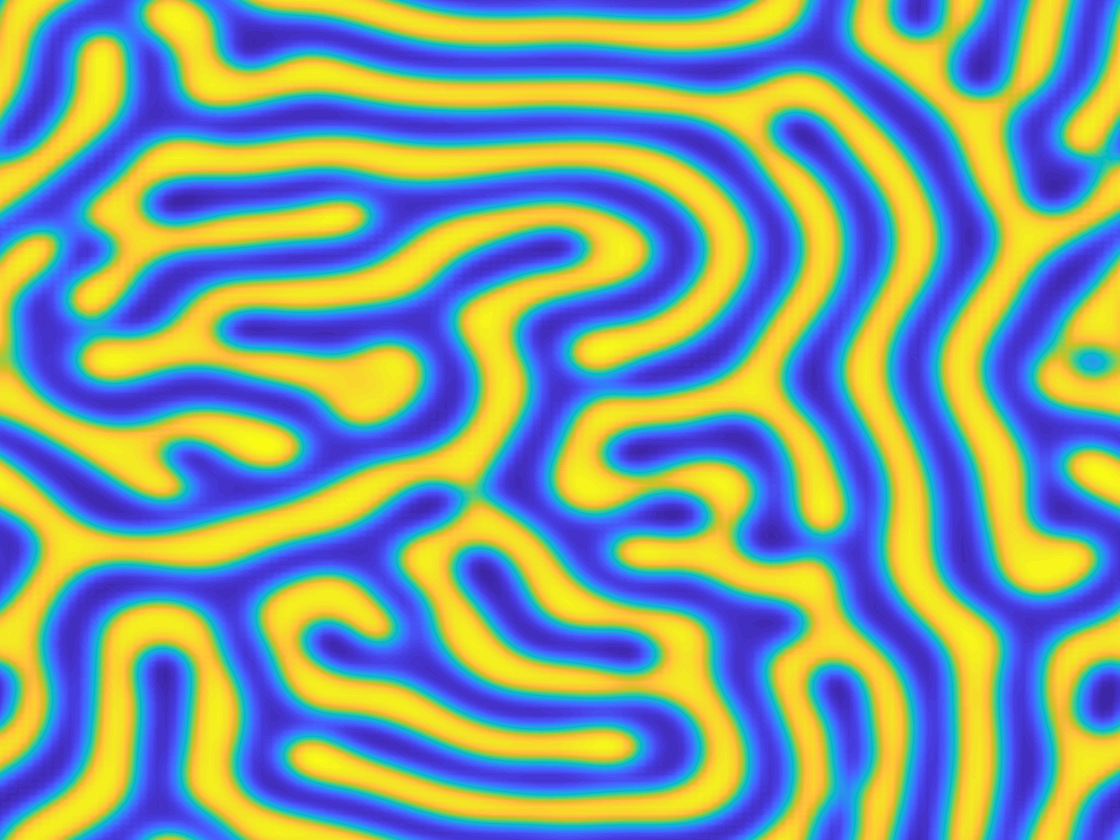}
		\includegraphics[width=1.40in]{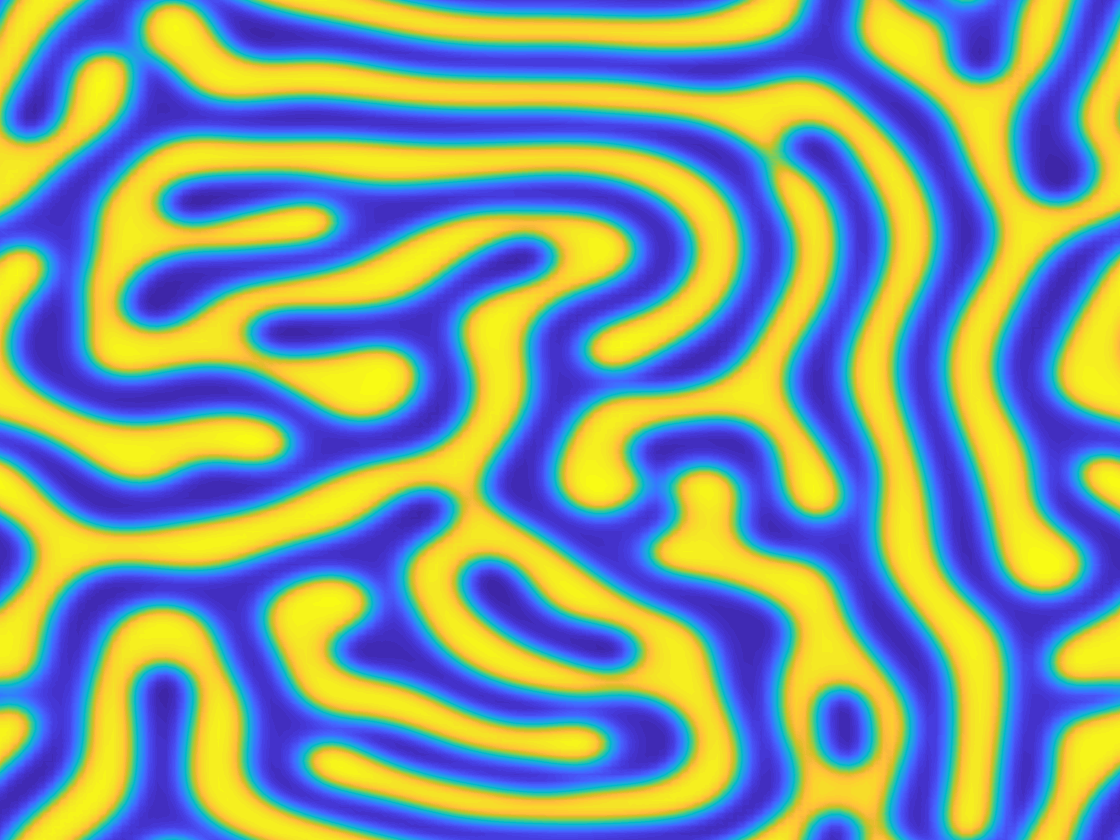}}
	\subfigure[The fractional order $\alpha=0.8$, at time $t=10,30,50,100$]{
		\includegraphics[width=1.40in]{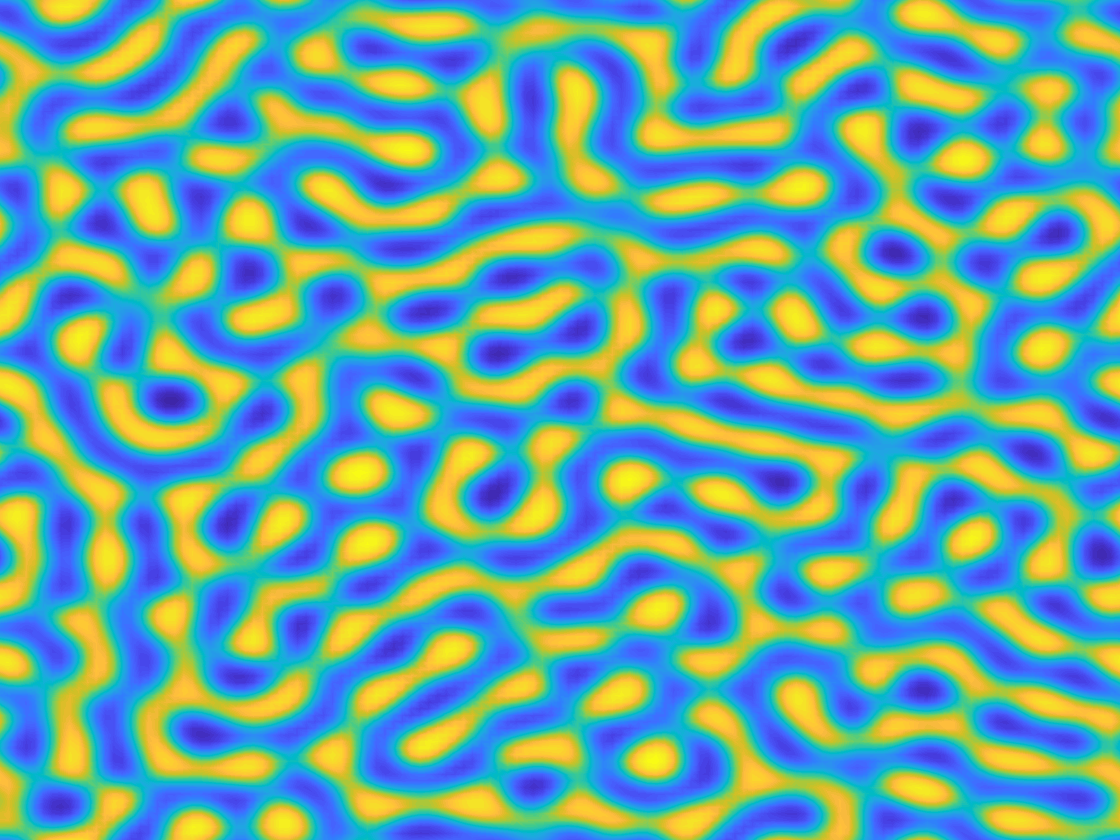}
		\includegraphics[width=1.40in]{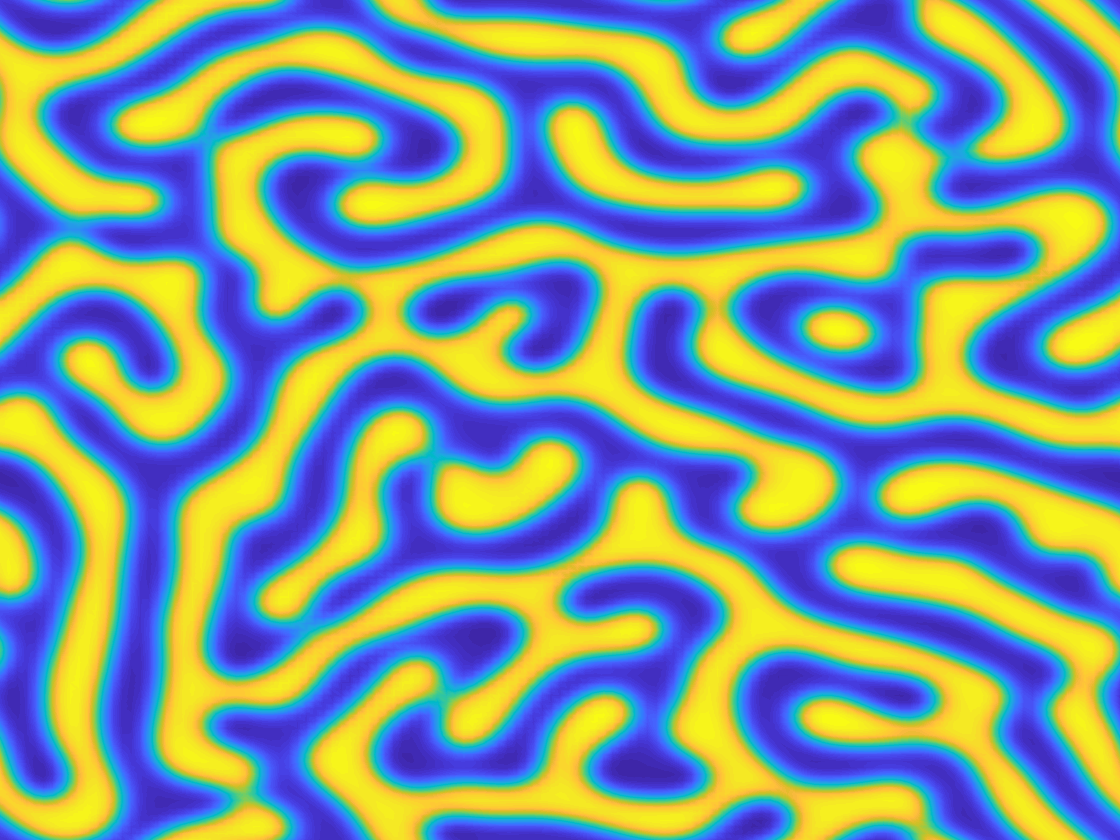}
		\includegraphics[width=1.40in]{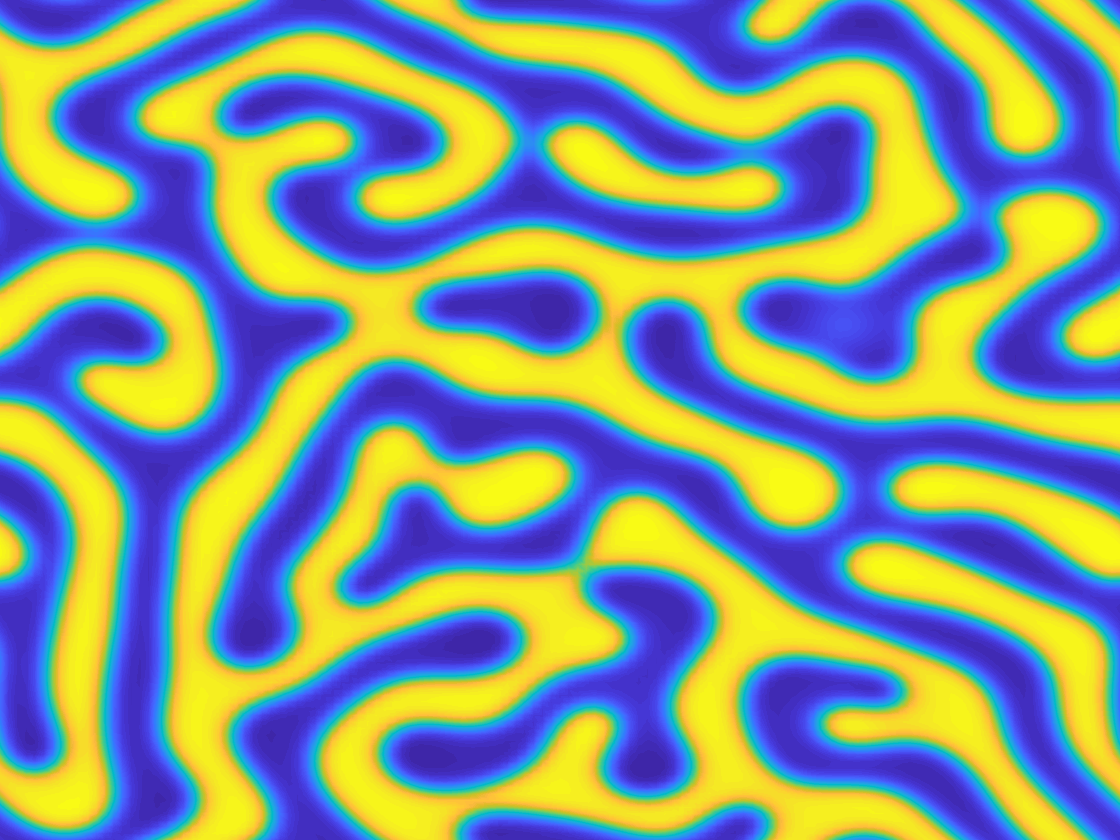}
		\includegraphics[width=1.40in]{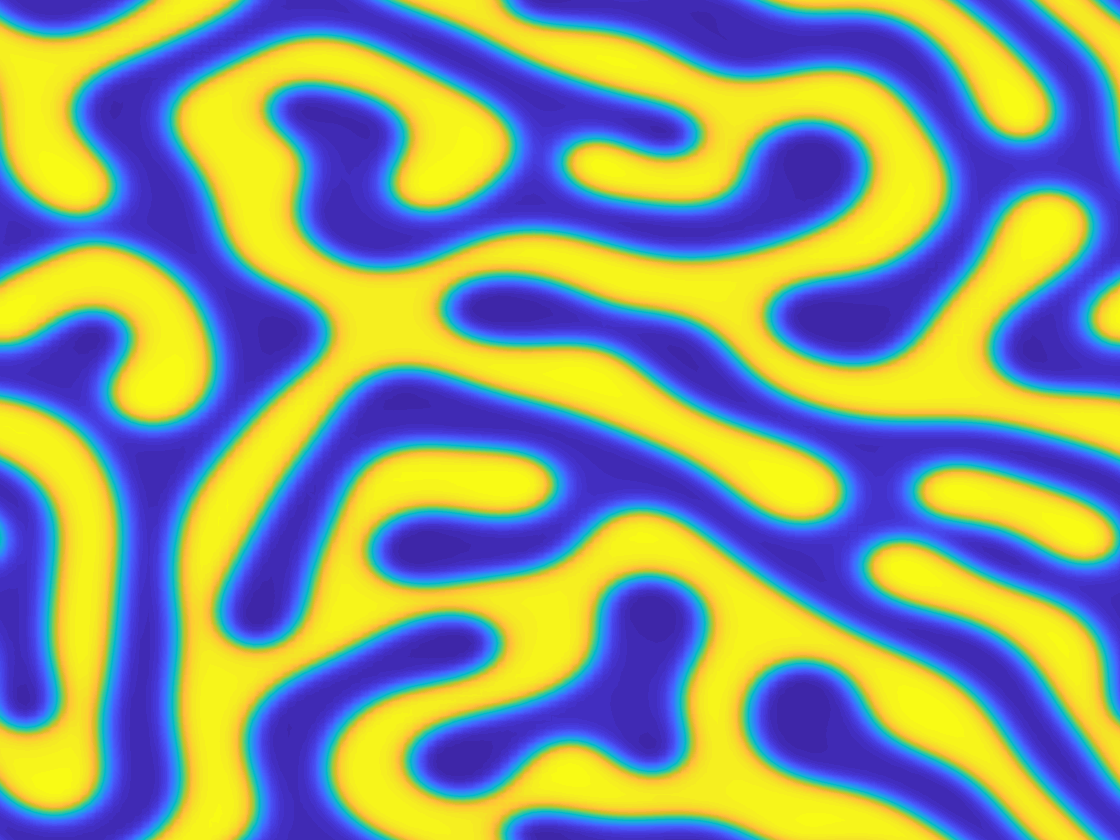}}
	\caption{Snapshots of dynamic coarsening processes for different fractional orders $\alpha$.}
	\label{figure:different alpha}
\end{figure}

To resolve the initial singularity, the graded time mesh $t_k=T_0(k/N_0)^{\gamma}$ with $T_0=0.01$,
$\gamma=2$ and $N_0=30$ is applied inside the initial period $[0,T_0]$.
Taking the predetermined maximum time-step $\tau_{\max}=10^{-1}$ and 
minimum step $\tau_{\min}=10^{-3}$, we use an adaptive time stepping strategy \cite{JiZhuLiao:2022TFCH}
in the remainder interval $[T_0,T]$,
\begin{align}\label{algorithm:adaptive time stepping}
	\tau_{ada}
	=\max\left\{\tau_{\min},
	\tau_{\max}/\Pi_{\phi}\right\}\quad\text{with}\quad \tau_{n+1}=\max\big\{\min\{\tau_{ada},R_*\tau_n\},r^*(\alpha)\tau_n\big\}
\end{align}
where $\Pi_{\phi}:=\sqrt{1+\eta\mynormb{\partial_\tau\phi^n}^2}$ and  $\eta>0$
is a user parameter.

Figure \ref{example:compar adap param} shows the numerical results by
using the adaptive time-stepping strategy \eqref{algorithm:adaptive time stepping}  
for  the fractional order $\alpha=0.5$ with three  parameters $\eta=10$, $10^2$ and $10^3$, respectively. 
In Figure \ref{example:compar adap param}, the uniform time 
step $\tau=5\times10^{-3}$ is used as a the reference solution.
Figure \ref{example:compar adap param}
plots the original energy curves and the associated time-steps. 
We choose the parameter $\eta=10^3$ for further simulations since it  
admits smaller time-steps in the fast varying region.

Figure \ref{example:FCH coarsen energy} depicts the numerical results by
using the adaptive time-stepping strategy \eqref{algorithm:adaptive time stepping} 
with the parameter $\eta=10^3$ until time $T=500$. The curves of original 
energy $E[\phi^n]$ in \eqref{cont:classical free energy} and the associated 
adaptive time steps are depicted in Figure \ref{example:FCH coarsen energy}. 
The adaptive time-stepping
approach seems to be efficient for our numerical simulations since the moment 
at which the energy decays rapidly can be captured with small time steps. 
As expected, the value of fractional order $\alpha$ significantly affects
the coarsening dynamics process and the energy curves always approach 
the integer-order one as the fractional order $\alpha\rightarrow1$.  
The phase profiles at four different times are listed in Figure \ref{figure:different alpha} for three  different fractional orders $\alpha=0.2,0.5,0.8$, respectivly. 
They are in accordance with the coarsening dynamics processes 
reported in \cite{JiZhuLiao:2022TFCH,LiuChengWangZhao:2018}.

\end{document}